\documentclass[reqno,oneside]{amsart}

\newcounter{heyheyCounter}[section]

\usepackage{amsthm}

\usepackage{
amssymb,amsmath,enumerate,stackrel,enumitem,
verbatim,
leftidx,bbm}
\usepackage{xcolor} \usepackage{xr-hyper}
\usepackage[all,cmtip]{xy}
\usepackage[utf8]{inputenc} \usepackage[T1]{fontenc}
\chardef\cprime"7E \usepackage{cases}

\usepackage[T1]{fontenc} 
\usepackage[pagebackref,hyperindex,linktocpage=true]{hyperref}
\hypersetup{
    colorlinks,
    linkcolor={red!50!black},
    citecolor={red!50!black},     urlcolor={red!50!black},     filecolor={red!50!black} }

\usepackage{theoremref}
\definecolor{labelkey}{rgb}{1,0,0}

\makeatletter
\@addtoreset{equation}{section}
\makeatother

\numberwithin{equation}{section}

\usepackage[pagebackref,hyperindex,linktocpage=true]{hyperref}
\hypersetup{
    colorlinks,
    linkcolor={red!50!black},
    citecolor={blue!50!black},
    urlcolor={blue!80!black}
}

\theoremstyle{definition}
\newtheorem{Defi}{Definition}[section] \newcommand{\defi}{\begin{Defi}} \newcommand{\xdefi}{\end{Defi}} \newtheorem{Cons}[Defi]{Construction} \newcommand{\const}{\begin{Cons}} \newcommand{\xconst}{\end{Cons}} \newtheorem{Hypo}[Defi]{Hypothesis} \newcommand{\hypo}{\begin{Hypo}} \newcommand{\xhypo}{\end{Hypo}} \newtheorem{Assu}[Defi]{Assumption} \newcommand{\assu}{\begin{Assu}} \newcommand{\xassu}{\end{Assu}} \newtheorem{DefiLemm}[Defi]{Definition and Lemma} \newcommand{\defilemm}{\begin{DefiLemm}} \newcommand{\xdefilemm}{\end{DefiLemm}} \newtheorem{Bsp}[Defi]{Example} \newcommand{\exam}{\begin{Bsp}} \newcommand{\xexam}{\end{Bsp}} \newtheorem{Syno}[Defi]{Synopsis} \newcommand{\syno}{\begin{Syno}} \newcommand{\xsyno}{\end{Syno}} \newtheorem{Bem}[Defi]{Remark} \newcommand{\rema}{\begin{Bem}} \newcommand{\xrema}{\end{Bem}} \newtheorem{Notation}[Defi]{Notation} \newcommand{\nota}{\begin{Notation}} \newcommand{\xnota}{\end{Notation}} 
\theoremstyle{plain}
\newtheorem{Theo}[Defi]{Theorem} \newcommand{\theo}{\begin{Theo}} \newcommand{\xtheo}{\end{Theo}} \newtheorem{Satz}[Defi]{Proposition} \newcommand{\prop}{\begin{Satz}} \newcommand{\xprop}{\end{Satz}} \newtheorem{Lemm}[Defi]{Lemma} \newcommand{\lemm}{\begin{Lemm}} \newcommand{\xlemm}{\end{Lemm}} \newtheorem{Coro}[Defi]{Corollary} \newcommand{\coro}{\begin{Coro}} \newcommand{\xcoro}{\end{Coro}}
\newtheorem{Ques}[Defi]{Question} \newcommand{\ques}{\begin{Ques}} \newcommand{\xques}{\end{Ques}}
\newtheorem{Conj}[Defi]{Conjecture} \newcommand{\conj}{\begin{Conj}} \newcommand{\xconj}{\end{Conj}}

\newcommand{\refit}[1]{(\ref{#1})}

\newcommand{\eqn}{\begin{equation}} \newcommand{\xeqn}{\end{equation}}
\newcommand{\eqnarr}{\begin{eqnarray*}} \newcommand{\xeqnarr}{\end{eqnarray*}}
\newcommand{\eqnarra}{\begin{eqnarray}} \newcommand{\xeqnarra}{\end{eqnarray}}

\newcommand{\pf}{\begin{proof}} \newcommand{\xpf}{\end{proof}}

\numberwithin{equation}{section}

\newcommand{\nc}{\newcommand}\nc{\StP}[1]{\cite[\href{http://stacks.math.columbia.edu/tag/#1}{Tag #1}]{StacksProject}}

\nc{\on}{\operatorname}
\nc{\aff}{{\on{aff}}}
\nc{\modi}{{\on{mod}}} \nc{\even}{{\on{even}}}
\nc{\odd}{{\on{odd}}}
\nc{\naive}{{\on{naive}}}
\nc{\hofib}{\on{hofib}}
\nc{\Bun}{\on{Bun}}
\nc{\ad}{{\on{ad}}}
\nc{\lft}{{\on{lft}}}

\nc{\Weil}{{\on{Weil}}} 
\nc{\FWeil}{{\on{FWeil}}} 

\nc{\str}{\on{-}}
\nc{\perf}{{\on{perf}}}
\nc{\Rel}{{\on{Pos}}}
\nc{\lan}{\langle}
\nc{\ran}{\rangle}
\nc{\tw}[1]{\langle #1 \rangle} 
\nc{\bbA}{{\mathbb A}} \nc{\bbB}{{\mathbb B}}
\nc{\bbC}{{\mathbb C}}
\nc{\bbD}{{\mathbb D}}
\nc{\bbE}{{\mathbb E}}
\nc{\bbF}{{\mathbb F}}
\nc{\bbG}{{\mathbb G}}
\nc{\bbH}{{\mathbb H}}
\nc{\bbI}{{\mathbb I}}
\nc{\bbJ}{{\mathbb J}}
\nc{\bbK}{{\mathbb K}}
\nc{\bbL}{{\mathbb L}}
\nc{\bbM}{{\mathbb M}}
\nc{\bbN}{{\N}} \nc{\bbO}{{\mathbb O}}
\nc{\bbP}{{\mathbb P}} \nc{\bbQ}{{\mathbb Q}} \nc{\bbR}{{\mathbb R}}
\nc{\bbS}{{\mathbb S}}
\nc{\bbT}{{\mathbb T}}
\nc{\bbU}{{\mathbb U}}
\nc{\bbV}{{\mathbb V}}
\nc{\bbW}{{\mathbb W}}
\nc{\bbX}{{\mathbb X}}
\nc{\bbY}{{\mathbb Y}}
\nc{\bbZ}{{\mathbb Z}}

\nc{\calA}{{\mathcal A}}
\nc{\calB}{{\mathcal B}}
\nc{\calC}{{\mathcal C}}
\nc{\calD}{{\mathcal D}}
\nc{\calE}{{\mathcal E}}
\nc{\calF}{{\mathcal F}}
\nc{\calG}{{\mathcal G}}
\nc{\calH}{{\mathcal H}}
\nc{\calI}{{\mathcal I}}
\nc{\calJ}{{\mathcal J}}
\nc{\calK}{{\mathcal K}}
\nc{\calL}{{\mathcal L}}
\nc{\calM}{{\mathcal M}}
\nc{\calN}{{\mathcal N}}
\nc{\calO}{{\mathcal O}}
\nc{\calP}{{\mathcal P}}
\nc{\calQ}{{\mathcal Q}}
\nc{\calR}{{\mathcal R}}
\nc{\calS}{{\mathcal S}}
\nc{\calT}{{\mathcal T}}
\nc{\calU}{{\mathcal U}}
\nc{\calV}{{\mathcal V}}
\nc{\calW}{{\mathcal W}}
\nc{\calX}{{\mathcal X}}
\nc{\calY}{{\mathcal Y}}
\nc{\calZ}{{\mathcal Z}}

\nc{\frakp}{{\mathfrak p}}
\nc{\frakm}{{\mathfrak m}}
\nc{\frakq}{{\mathfrak q}}

\nc{\unit}{1}

\nc{\Sht}{{\on{Sht}}}
\nc{\Frob}{{\on{Frob}}}
\nc{\Hecke}{{\on{Hecke}}}
\nc{\inv}{{\on{inv}}}
\nc{\Conv}{{\on{Conv}}}
\nc{\triv}{{\on{triv}}}
\nc{\Isom}{{\on{Isom}}}

\nc{\scrB}{{\mathscr{B}}}
\nc{\scrA}{{\mathscr{A}}}
\nc{\bbf}{{\mathbf{f}}}
\nc{\bba}{{\mathbf{a}}}
\nc{\rig}{{\mathrm rig}}

\nc{\al}{\alpha}
\nc{\be}{\beta}
\nc{\ga}{\gamma}
\nc{\la}{\lambda}
\nc{\qcqs}{{\on{qcqs}}}
\nc{\Bmu}{{\boldsymbol \mu}}

\nc{\pot}[1]{ [\hspace{-0,5mm}[ {#1} ]\hspace{-0,5mm}] }
\nc{\rpot}[1]{ (\hspace{-0,7mm}( {#1} )\hspace{-0,7mm}) }

\nc{\defined}{\hspace{0.1cm}\stackrel{\text{\tiny \rm def}}{=}\hspace{0.1cm}}
\nc{\co}{\colon}
\nc{\specto}{{\leadsto}}

 \newcommand{\category}[1]{\mathrm{#1}}

    \newcommand{\Cat}{\category{Cat}}   \newcommand{\Ho}{\category{Ho}}           \newcommand{\HypShv}{\category{HypShv}}  \newcommand{\Fun}{\category{Fun}}   \newcommand{\PreStk}{\category{PreStk}}    \newcommand{\Ani}{\category{Ani}}      \newcommand{\Spt}{\category{Spt}}      \newcommand{\Sm}{\category{Sm}}  \newcommand{\Corr}{\category{Corr}} \newcommand{\Sch}{\category{Sch}}      \renewcommand{\Pr}{\category{Pr}}
\newcommand{\PrL}{\Pr^\category{L}} \newcommand{\PrgmL}{\PrL_\mathrm{gm}}    \newcommand{\st}{{\category{st}}}   \newcommand{\PrSt}{{\Pr_{\category{st}}^{\category{L}}}} \newcommand{\PrLSt}{{\Pr_{\st}^{\category{L}}}}      \newcommand{\PrStOm}{{\Pr_{\omega,\st}^{\category{L}}}} \newcommand{\PrLomega}{\Pr_\omega^{\category L}}
\newcommand{\PrLomegast}{\Pr_{\omega, \st}^\category{L}}
\newcommand{\PrLomegaZ}{\Pr_{\omega, \Z}^\category{L}}
\nc{\PrLStomega}{\PrLomegast}
\nc{\PrLstomega}{\PrLomegast}

\nc{\idem}{\category{idem}}

 \newcommand{\dbl}{\category{dbl}}  \newcommand{\Ind}{\operatorname{Ind}}       \newcommand{\Mod}{\category{Mod}}                 \newcommand{\W}{\mathrm {W}}  \newcommand{\Dcons}{{\category{D}_{\textup{cons}}}}

\newcommand{\DAf}{{\DA^\flat}}

         \def\Gm{\mathbf {G}_\mathrm m}         \def\red{\mathrm{red}} \def\bk{{\bar k}}   \def\ft{\mathrm{ft}}     
 
   \newcommand{\colim}{\operatornamewithlimits{colim}}  \newcommand{\cofib}{\operatorname{cofib}}       \def\id{{\rm id}} \def\ev{{\operatorname {ev}}} \def\pr{{\rm pr}} \def\opp{{\rm op}}  \def\To#1#2{\mathop{\count0=#1 \loop\ifnum\count0>0 \smash-\mkern-7mu \advance\count0 -1 \repeat \mathord\rightarrow}\limits^{#2}}    \def\Map{\mathop{\rm Map}\nolimits}  \def\Char{\mathop{\rm char}\nolimits}    \def\Char{\mathop{\rm char}} \def\Hom{\mathop{\rm Hom}\nolimits}   \def\coev{\mathop{\rm coev}\nolimits}      
     \def\Sht{\mathop{\rm Sht}\nolimits}  \def\CAlg{\mathop{\rm CAlg}\nolimits}    \def\IHom{\underline{\Hom}}  \def\End{\mathop{\rm End}\nolimits}     \def\Map{\mathop{\rm Map}\nolimits}          \def\et{\mathrm{\acute et}}      
\definecolor{hellgrau}{RGB}{200,200,200} \definecolor{dunkelgrau}{RGB}{160,160,160} \definecolor{hellblau}{RGB}{194, 215, 249} \definecolor{dunkelblau}{RGB}{68, 128, 226} 
\def\Z{{\mathbb{Z}}} \def\F{{\mathbb F}}  \def\Fpq{{\ol \F_p}} \def\Fq{{\F_q}}   \def\N{{\mathbb N}} \def\Q{{\mathbb Q}}   \def\Zhatp{{\hat\Z^p}}
  \def\Zl{{\Z_\ell}}   \def\A{{\mathbb A}}  \def\proet{{\textup{proét}}}  \def\Gm{{\mathbb {G}_\mathrm m}}

\def\H{{\rm H}}     \def\SH{\category{SH}} \def\nc{\mathrm{nc}}    \def\gm{\mathrm{gm}} \newcommand{\cons}{\mathrm{cons}}
\newcommand{\indbl}{\category{ind\textrm{-}dbl}}

\def\Dad{\widehat {\category{D}}}     
\def\DMflat{\DA^\flat}
\def\DAf{\DMflat}

\def\eff{\mathrm{eff}} \def\DA{\category{DA}}         \def\ii{$\infty$}

  \def\RG{\R \Gamma}    

\def\Spec{\mathop{\rm Spec}}             \newcommand{\comp}{\omega} \newcommand{\at}{\mathrm{at}}           \newcommand{\D}{\category{D}} \newcommand{\tD}{\widetilde{\category{D}}} \newcommand{\DW}{\D^W}    
   
\def\R{{\rm R}}   \def\sbuildrel#1\over#2{\mathrel{\smash{\mathop{\kern0pt #2}\limits^{#1}}}}

\def\dual{\vee}

\let\x\times
\let\ol\overline

\renewcommand{\t}{\otimes}

\newcommand{\bx}{\boxtimes}

\renewcommand{\r}{\rightarrow}
\newcommand{\lr}{\longrightarrow}

\begin{document}

\title[]{Categorical K\"unneth formulas \\ for cohomological motives}
\author[]{Timo Richarz, Jakob Scholbach}

 \thanks{
The first named author T.R.~is funded by the European Research Council (ERC) under Horizon Europe (grant agreement nº 101040935), by the Deutsche Forschungsgemeinschaft (DFG, German Research Foundation) TRR 326 \textit{Geometry and Arithmetic of Uniformized Structures}, project number 444845124 and the LOEWE professorship in Algebra, project number LOEWE/4b//519/05/01.002(0004)/87.
The second named author J.S. acknowledges supported by the European Union -- Project 20222B24AY (subject area: PE -
Physical Sciences and Engineering) ``The arithmetic of motives and L-functions'', and logistical support by the Max-Planck-Institute for Mathematics in Bonn.
Both authors acknowledge support by the MSCA doctoral network ReMoLD (project number 101168795) under Horizon Europe.
}

\address{Fachbereich Mathematik, TU Darmstadt, Schlossgartenstrasse 7, 64289 Darmstadt, Germany}
\email{richarz@mathematik.uni-darmstadt.de}

\address{Università degli Studi di Padova, Via Trieste 63, Padova, Italia}

\maketitle

\begin{abstract} 
The manuscript at hand systematically studies K\"unneth formulas at a categorical level.
We give criteria for an abstract six functor formalism to satisfy the categorical K\"unneth formula, and use this to formulate conjectures for categories of étale motives.
As supporting evidence for these conjectures, we prove categorical Künneth formulas for adic sheaves and for cohomological motives, i.e., étale motives modulo the kernel of the adic realization.
\end{abstract}

\tableofcontents

\section{Introduction}
Given two geometric objects $X$ and $Y$ it is natural to investigate how invariants of their product $X \x Y$ are related to the ones on the factors. 
For example, for two complex algebraic varieties $X$ and $Y$ the \textit{K\"unneth formula} gives a canonical isomorphism of their Betti cohomology groups
\begin{equation}
\H^n(X \x Y, \Q) = \bigoplus_{i+j=n} \H^i(X, \Q) \t_\Q \H^j(Y, \Q).
\label{Künneth cohomology}
\end{equation}
In other words, every cohomology class on the product can be written as a linear combination of cohomology classes on the factors. 
Similar K\"unneth formulas also hold true for (algebraic) de Rham cohomology of smooth varieties over fields of characteristic $0$ and $\ell$-adic cohomology over algebraically closed fields of characteristic different from $\ell$.

Motivic cohomology, however, fails to satisfy a Künneth formula as in \eqref{Künneth cohomology}.
Indeed, for a fixed algebraic variety $X$, Totaro has shown that the K\"unneth formula holds for its rationalized higher Chow groups and all smooth proper varieties $Y$ if and only if $X$ is mixed Tate \cite[Theorems~4.1, 7.2]{Totaro:Motive}.
The only smooth proper connected curve $X$ over an algebraically closed field that satisfies this is $X = \mathbb P^1$.

We propose to sidestep this failure by considering a Künneth formula at the level of the category of motives, as opposed to the level of individual cohomology groups (Conjectures~\ref{conjecture kuenneth intro}~and~\ref{conjecture kuenneth weil intro}).
This approach is inspired by categorical Künneth formulas for, say, derived categories of quasi-coherent sheaves or $D$-modules \cite{Lurie:HA, BenZviFrancisNadler:Integral, DrinfeldGaitsgory:Some}.
From this perspective, the present manuscript investigates the extent to which these formulas remain valid for categories of motives. 

Our approach builds on a systematic study of the categorical K\"unneth formula for an abstract $6$-functor formalism (Section~\ref{section:abstract-kuenneth}).  
The main results establish such formulas for ind-dualizable cohomological motives (\thref{main thm intro}) and their variants for Weil objects (\thref{kuenneth DAf weil intro}). 

\subsection{A conjectural Künneth formula for étale motives}
Let $k$ be a field of finite Galois cohomological dimension, e.g.,  a separably closed field, a finite field or an unordered number field. 
For any algebraic $k$-variety $X$, we denote by $\DA(X)$ the category of étale motives from \cite{Ayoub:Realisation,CisinskiDeglise:Etale}.
Briefly, this category is constructed from $\mathbb A^1$-invariant sheaves on the smooth-étale site of $X$ by inverting the Tate twist.
After rationalization, Chow groups may be computed as homomorphism groups in $\DA(X)$ when $X$ is regular.  
Étale motives come equipped with a $6$-functor formalism. 

Taking advantage of its natural \ii-categorical enhancement, $\DA(X)$ is a module over the presentably symmetric monoidal category $\DA(k)$ (Section~\ref{section:algebraic-motives}).
We can thus form the relative Lurie tensor product (Section~\ref{section:infinity-cats-notation}), and consider the \textit{exterior product functor}
\begin{equation}\label{kuenneth functor intro}
\DA(X) \t_{\DA(k)} \DA(Y) \r \DA(X \x_k Y)
\end{equation}
induced from the map $(M,N)\mapsto M\bx N$ for objects $M\in \DA(X), N\in \DA(Y)$.
Here, we denote by
$$ M\bx N:=\pr_X^*M\t \pr_Y^*N\in \DA(X \x_k Y)$$
the exterior product where $\pr_X$ and $\pr_Y$ denote the projections from $X\x_k Y$ to $X$ and $Y$, respectively. 

\prop[{\thref{proposition:SH-fully-faithful}}]
The functor \eqref{kuenneth functor intro} is fully faithful.
\xprop

This is a consequence of the Künneth formula for $M\bx N$ due to Jin--Yang \cite{JinYang:Kuenneth} and the computation of enriched maps in Lurie tensor products (\thref{proposition:infinity-cats-tensor-product}). 
The latter uses compact, more precisely, atomic generation of these categories which is ensured by the Galois cohomological assumption on $k$.
A similar statement (\thref{fully faithful SH}) holds for the stable $\A^1$-homotopy category $\SH$ in place of $\DA$ (for general fields $k$), and also for categories of adic sheaves $\Dad$ and of cohomological motives $\DAf$ considered in Section~\ref{coho mot intro}.

The \textit{categorical K\"unneth formula} asks, in addition to full faithfulness, for the essential surjectivity of the exterior product functor, i.e., whether any étale motive on the product variety $X\x_k Y$ can be decomposed as a colimit of motives coming from the factors $X$ and $Y$:

\conj[{\thref{conjecture kuenneth}}]\thlabel{conjecture kuenneth intro}
Assume $Y$ is smooth over $k$ and, in addition, one of the following:
\begin{enumerate}
    \item \label{char0 intro} $\Char k=0$, or
    \item \label{Y proper intro}$Y$ is proper over $k$.
\end{enumerate}
Then, the exterior product functor \eqref{kuenneth functor intro} induces an equivalence:
$$\DA(X)^\indbl \t_{\DA(k)} \DA(Y)^\indbl \stackrel {\cong\, ?} \lr \DA(X \x_k Y)^\indbl,$$
where the superscript $\indbl$ (=ind-dualizable) denotes the subcategories generated under colimits by dualizable objects.
\xconj

The restriction to ind-dualizable objects is analogous to the passage from constructible sheaves to local systems.
It seems necessary since objects with diagonal support generally do not lie in the essential image of \eqref{kuenneth functor intro}. 

Our approach towards the above conjecture is based on a construction due to Gaitsgory--Kazhdan--Rozenblyum--Varshavsky \cite[A.5.7]{GaitsgoryEtAl:Toy}. 
We prove that the above conjecture for $X$ and $Y$ will be implied by the following two conditions (\thref{purity remark}):
\begin{enumerate}
    \item[(a)]
    For any \emph{closed} point $x \in X$, the pullback functor
    $$
(x\x \id_Y)^*\colon \DA(X\x_k Y)^\indbl\r \DA(x\x_k Y) 
$$
is conservative if $X$ is connected. 
    \item[(b)] The pushforward $\pr_{X,*}\colon \DA(X \x_k Y) \r \DA(X)$ along the projection $X \x_k Y \r X$ preserves dualizable objects.
\end{enumerate}

Colloquially speaking, (a) is asking for $\DA(-)^\indbl$ to satisfy a Tannakian formalism with pullbacks to closed points serving as fiber functors.
This conservativity is implied by a strong form of the conservativity of the motivic nearby cycles functor 
conjectured by Ayoub \cite[Conjecture~5.1]{Ayoub:Motivic}.
For étale \emph{torsion} motives (equivalently, étale torsion sheaves), it holds by classical theory.

Condition (b) holds if $Y$ is smooth and proper over $k$, which explains the appearance of Condition \eqref{Y proper intro} in \thref{conjecture kuenneth intro}.
In sync with Condition \eqref{char0 intro} there, we prove that for $Y$ smooth but not necessarily proper, and $\Char k = 0$, dualizable étale \emph{torsion} motives are preserved under pushforward (\thref{pushforward char 0}).
For $X = Y = \A^1_k$ and $\Char k > 0$, Condition (b) is not satisfied (\thref{A2}).
This problem can be fixed by considering Weil objects (Section~\ref{weil intro}).  

Conditions \eqref{char0 intro} and \eqref{Y proper intro} in \thref{conjecture kuenneth intro} match the conditions in \cite{SGA1} under which the map of étale fundamental groups
$$\pi_1(X \x_k Y) \r \pi_1 (X) \x \pi_1(Y)$$
is an isomorphism, with $k$ algebraically closed.

\subsection{The Künneth formula for cohomological motives and adic sheaves}\label{coho mot intro}
As supporting evidence for \thref{conjecture kuenneth intro},
we prove the categorical Künneth formula for adic sheaves $\Dad$ and cohomological motives $\DAf$ in place of $\DA$.
For any $k$-variety $X$, denote by $\Dad(X)$ the profinite completion of $\DA(X)$. 
As explained in Section~\ref{section:adic-sheaves}, this agrees with more classical versions of categories of adic sheaves by the rigidity theorem. 
(Further completion at some prime $\ell\neq \Char k$ recovers the category of $\ell$-adic sheaves.)

The category comes equipped with an \textit{adic realization functor}
$$\rho\colon \DA(X)\r \Dad(X),$$
given on compact objects by $\rho(M)=\lim_{n\in \Z_{\geq 1}}M/n$.
The following is analogous to a construction due to Ayoub \cite{Ayoub:Motives}, who used the Betti realization instead.

\defi
The category of \textit{cohomological motives} $\DAf(X)$ is the Verdier quotient of $\DA(X)$ by the subcategory of objects $M\in \DA(X)$ with $\rho(M)\simeq 0$.
\xdefi

By construction, the category of cohomological étale motives $\DAf(X)$ fits into a diagram
$$\xymatrix{
\DA(X) \ar[rr]^{\rho} \ar[dr]_{(-)^\flat} & & \Dad(X) \\
& \DAf(X) \ar[ur]_{\rho^\flat},}$$
where $(-)^\flat$ is a localization and $\rho^\flat$ is conservative. 
Thus, some $M \in \DA(X)$ is zero in $\DAf(X)$ precisely if and only if $\rho(M) \simeq 0$.

\theo[{Section~\ref{section:six-functors}}]\thlabel{six functor DAf intro}
The assignment $X \mapsto \DAf(X)$ is part of a $6$-functor formalism satisfying the properties familiar for $\DA$ or $\Dad$, including excision and relative purity.
In addition, each category $\DAf(X)$ is dualizable as a module over $\DAf(k)$.
\xtheo

In comparison to these good abstract properties, we have little control on the kernel of $\rho$ and thus on the morphism groups in $\DAf(X)$. 
Even though $\rho$ preserves compact objects, we only know that $\DAf(X)$ is $\omega_1$-compactly generated. 
In particular, we do not know whether $\DAf(k)$ is compactly generated, nor whether it is dualizable as a module over spectra.
However, the above-mentioned dualizability of $\DAf(X)$ over $\DAf(k)$, which we prove using Ramzi's work on dualizable categories \cite{Ramzi:DualizableCategories}, can be seen as a replacement for an absolute dualizability property.

Of course, whenever $\rho$ is conservative, as one might conjecture when $k$ is of finite transcendence degree over its prime field \cite{Ayoub:Motives}, we will have $\DA(X) = \DAf(X)$.
For an application towards \thref{conjecture kuenneth intro}, a weaker conjecture suffices, namely conservativity of $\rho$ when restricted to ind-dualizable and compact objects \cite{ayoub:weilcohomologytheoriesmotivic}.

\theo[{\thref{Künneth DMf indlis}}]\thlabel{main thm intro}
Under either of the assumptions in \thref{conjecture kuenneth intro}, ind-dualizable cohomological motives satisfy the categorical Künneth formula, i.e., the exterior product functor induces an equivalence:
$$\DAf(X)^\indbl \t_{\DAf(k)} \DAf(Y)^\indbl \stackrel \cong \lr \DAf(X \x_k Y)^\indbl.$$
\xtheo

A similar equivalence holds for ind-dualizable adic sheaves $\Dad(-)^\indbl$ in place of $\DA(-)^\indbl$ (\thref{Künneth Dad indlis}).

\subsection{The case of Weil motives}\label{weil intro}
For the remainder of the introduction, let $k=\Fq$ be a finite field of characteristic $p>0$ and cardinality $q$, and fix an algebraic closure $\bar k$.
For any $k$-variety $X$ denote by $\ol X:=X\x_k\bar k$ its base change equipped with the $\bar k$-endomorphism $\varphi_X:=\Frob_X\x \id_{\bar k}$ where $\Frob_X$ is the $q$-Frobenius on $X$.
 
We define the category of \textit{étale Weil motives} $\DA(X^W)$ to be the ind-completion of the categorical fixed points of the automorphism $\varphi_X^*$ on the compact objects $\DA(\ol X)^\omega$. 
Compact objects in this category are pairs of objects $M \in \DA(\ol X)$ with an isomorphism $\varphi_X^* M \cong M$. 
The category $\DA(X^W \x_{\bar k} Y^W)$ is defined similarly, using two commuting isomorphisms for $\varphi_X^*$ and $\varphi_Y^*$.
In Section~\ref{sect--Künneth Weil motives}, we show that the assignment $X\mapsto \DA(X^W)$ upgrades to a $6$-functor formalism, linearly over $\DA(\bar k)$.

\conj[{\thref{conjecture kuenneth weil}}]\thlabel{conjecture kuenneth weil intro}
For any $k$-varieties $X$ and $Y$, the exterior product functor induces an equivalence on the categories of étale Weil motives:
$$\DA(X^\W) \t_{\DA(\bar k)} \DA(Y^W) \stackrel {\cong\, ?} \lr \DA(X^W \x_{\bar k} Y^W).$$
\xconj

Again, our evidence is a result for cohomological Weil motives, constructed in analogy to $\DAf(X)$ by passing to the Verdier quotient by the kernel of the adic realization.

\theo[{\thref{Künneth DMf Weil}}]\thlabel{kuenneth DAf weil intro}
Cohomological Weil motives satisfy the categorical K\"unneth formula, i.e., for any $k$-varieties $X$ and $Y$ the exterior product functor induces an equivalence: 
$$\DAf(X^W) \t_{\DAf(\bar k)} \DAf(Y^W) \stackrel \cong \lr \DAf(X^W \x_{\bar k} Y^W).$$
\xtheo

This result, which also holds for categories of adic sheaves $\Dad$ in place of $\DAf$ (cf.~also \cite{HemoRicharzScholbach:Kuenneth}), potentially allows to decompose the intersection motive of the moduli stack of shtukas, constructed motivically in \cite{RicharzScholbach:Intersection}, as a cohomological motive.

\medskip
\noindent\textbf{Acknowledgements.} 
The manuscript at hand profited from discussions with Joseph Ayoub, Chirantan Chowdhury, Rızacan Çiloğlu, Denis-Charles Cisinski, Dennis Gaitsgory, Marc Hoyois, Markus Land, Thomas Nikolaus, Maxime Ramzi, Peter Scholze, Fabio Tanania and Can Yaylali.

\section{Preliminaries on \ii-categories}\label{section:infinity-cats}

The categorical Künneth formula for cohomological motives in Theorems~\ref{main thm intro}~and~\ref{kuenneth DAf weil intro} involves a tensor product over a symmetric monoidal category generated by dualizable, but not necessarily compact objects.
In this section, we provide a tool to compute enriched mapping objects in such tensor products (\thref{proposition:infinity-cats-tensor-product}) by drawing on Ramzi's work on dualizable categories \cite{Ramzi:DualizableCategories}. 
We also study how to enforce a functor to be conservative (Section~\ref{sect--enforcing conservativity}).

\subsection{Notation}\label{section:infinity-cats-notation}

The category of presentable \ii-categories with colimit-preserving (or, continuous) functors $\PrL$ is closed symmetric monoidal with respect to the tensor product $\t$ defined in \cite[Section~4.8]{Lurie:HA}.

Let $\CAlg(\PrL)$ be the category of commutative algebra objects in $\PrL$, i.e., presentably symmetric monoidal \ii-categories with colimit-preserving, symmetric monoidal functors. 
For $A \in \CAlg(\PrL)$, we denote by $\Mod_A(\PrL)$ the category of $A$-modules in $\PrL$.
It carries a symmetric monoidal structure $\t_A$ with unit $A$. 
The tensor product $C\t_A D\in \Mod_A(\PrL)$ is computed from the tensor product in $\PrL$ using a two-sided bar construction \cite[Theorem~4.4.2.8]{Lurie:HA}.
The forgetful functor $\Mod_A(\PrL)\to \PrL$ is lax symmetric monoidal, and preserves limits and colimits \cite[Corollaries 4.2.3.3, 4.2.3.5]{Lurie:HA}.

Recall from \cite[§4]{Lurie:HA} that an object $x$ in a symmetric monoidal \ii-category is called \emph{dualizable} if there is another object $x^\dual$ and so-called \emph{coevaluation} and \emph{evaluation} maps 
$$\coev : 1 \r x \t x^\dual,\;\; \ev : x \t x^\dual \r 1$$
satisfying the triangle identities.

\subsection{Enriched maps}
For an \ii-category $C$, we denote by
\begin{equation}
    \Map_C \colon C^\opp \x C \r \Ani \label{mapping functor}
\end{equation}
the mapping functor, i.e., $ \Map_C(c,c')$ is the anima of maps (or, mapping space) for objects $c,c'\in C$. 
For $A\in \CAlg(\PrL)$ and $C\in \Mod_A(\PrL)$, we have the $A$-enriched mapping functor
\eqn\label{enriched mapping functor}
\Map_{C/A}\colon C^\opp \x C \r A.
\xeqn
By definition, $\Map_A(-,\Map_{C/A}(c,c')) = \Map_C(-\t c,c')$ as functors $A^\opp \r \Ani$.
Note that $\Map_{C/A}(c,c')$ exists by the adjoint functor theorem and that \eqref{enriched mapping functor} constitutes a functor by \cite[Proposition~5.1.10]{Land:Book}.

\exam\thlabel{example:enriched-hom-spectra}
The enriched mapping functor interpolates between the following two extreme cases:
\begin{enumerate}
    \item For anima $A = \Ani$, we have $\Map_{C/\Ani} = \Map_C$ as in \eqref{mapping functor}. For spectra $A = \Spt$, we get the mapping spectrum, so $\pi_0 \Map_{C/\Spt}$ is the set of morphisms in the homotopy category $\Ho(C)$, which is a triangulated category.
    \item If $C=A$, then $\Map_{A/A} = \IHom_A$, i.e., the inner hom in $A$.\end{enumerate}
\xexam

\lemm\thlabel{lemma:fully-faithful-enriched}
A functor $f\colon C\to D$ in $\Mod_A(\PrL)$ is fully faithful if and only if the induced map on $A$-enriched maps
\eqn\label{lemma:fully-faithful-enriched:eq}
\Map_{C/A}(c,c')\to \Map_{D/A}(f(c),f(c'))
\xeqn
is an isomorphism in $A$ for all $c,c'\in C$.
\xlemm
\pf
The map \eqref{lemma:fully-faithful-enriched:eq} exists by $A$-linearity of $f$.
If \eqref{lemma:fully-faithful-enriched:eq}  is an isomorphism (in $A$), we can apply $\Map_A(1_A, -)$ and recover the mapping spaces, so $f$ is fully faithful. 
The converse follows from applying the Yoneda embedding $A^\opp \to \Fun(A,\Ani), a\mapsto \Map_A(a,-)$ to the mapping anima and using the $A$-linearity of $f$.
\xpf

\subsection{Internal left adjoints}\label{internal left adjoint section}
Throughout this subsection, let $A \in \CAlg(\PrL)$ and $C \in \Mod_A(\PrL)$.
In particular, the action functor $A \x C \r C$, $(a, c) \mapsto a \t c$ preserves colimits in both variables.

\defi\thlabel{internally left adjoint}
A functor $f\colon C\to D$ in $\Mod_A(\PrL)$ is \emph{internally left adjoint} if its right adjoint $f^R$ (which a priori is only lax $A$-linear and need not preserve colimits) lies in $\Mod_A(\PrL)$ as well.
\xdefi

The following notion is used in work of Ben Moshe \cite[Section~5]{Ben-Moshe:Yoneda} and Ramzi \cite[Section~1.2]{Ramzi:DualizableCategories}.
We will use it to control enriched mapping objects in tensor products (\thref{proposition:infinity-cats-tensor-product}).

\defi\thlabel{atomic objects}
An object $c\in C$ is called \textit{$A$-atomic} (or simply \textit{atomic}) if the functor $A \r C, a \mapsto a \t c$ is internally left adjoint, i.e., if $\Map_{C/A}(c,-): C\r A$ is $A$-linear and preserves colimits.
Denote by $C^{\textup{at}/A}\subset C$ (or simply $C^{\textup{at}}$) the full subcategory of atomic objects, which is a small category \cite[Corollary~1.26]{Ramzi:DualizableCategories}.

In addition, $C$ is \textit{atomically generated by $C_0\subset C^{\textup{at}}$} if the smallest $A$-submodule of $C$ that is closed under colimits and contains $C_0$ is $C$ itself.
\xdefi

\exam\thlabel{atomic example}
The notion interpolates between the following cases:
\begin{enumerate}
	\item \label{atomic example spectra} For spectra $A = \Spt$, the atomic objects in $C$ are exactly the compact objects \cite[Example~1.24]{Ramzi:DualizableCategories} (see \thref{atomic vs compact} for a generalization). 
    	\item \label{atomic example dualizable} If $C=A$, the atomic objects in $A$ are exactly the dualizable objects \cite[Example~1.23]{Ramzi:DualizableCategories}.
\end{enumerate}
\xexam

\lemm\thlabel{internal left adjoint lemma}
Let $f\colon C\r D$ be a functor in $\Mod_A(\PrL)$.
If $f$ is internally left adjoint, then $f$ preserves atomic objects. 
Conversely, if $C$ is atomically generated by $C_0\subset C^{\textup{at}}$ and $f(C_0)\subset D^{\textup{at}}$, then $f$ is internally left adjoint. 
In this case, $f$ is fully faithful if and only if for all objects $c,c'\in C_0$:
\eqn\label{lemma:fully-faithful-relative-compact:eq}
\Map_{C/A}(c,c')\overset{\cong}{\lr} \Map_{D/A}(f(c),f(c'))
\xeqn
\xlemm
\pf
The characterization of internal left adjoints is \cite[Corollary~1.31]{Ramzi:DualizableCategories}.
Further, the unit of the adjunction $\id_C\r f^R\circ f$ is an equivalence if its restriction to $C_0$ is so.
Then \eqref{lemma:fully-faithful-relative-compact:eq} is an equivalence for all $c,c'\in C$, which is equivalent to the full faithfulness of $f$ by \thref{lemma:fully-faithful-enriched}.
\xpf

Under additional assumptions on $A$ and $C$ there is a close relation between atomic and compact objects:

\lemm
\thlabel{atomic vs compact}
Suppose $A$ is stable and generated under colimits by dualizable objects. 
Then, the following hold:
\begin{enumerate}
\item \label{atomic vs compact 1} An object $c\in C$ is atomic if and only if $\Map_{C/A}(c,-) \colon C\r A$ preserves filtered colimits (equivalently, direct sums).
\item \label{atomic vs compact 2} If $A$ is compactly generated and the action $A\x C\r C$ preserves compact objects, then every compact object in $C$ is atomic.
\end{enumerate}
In particular, if $A$ is locally rigid \cite[Definition~D.7.4.1]{Lurie:SAG} and $A\x C\r C$ preserves compact objects, then $C^{\textup{at}}=C^\omega$, i.e., atomic objects are exactly the compact ones.
\xlemm
\pf
Locally rigid categories are compactly generated and the monoidal unit $1\in A$ is compact.
The latter implies that atomic objects are compact. 
So, the final statement follows from Item \eqref{atomic vs compact 2}. 

\textit{Item \eqref{atomic vs compact 1}, ``if'':} Since $f=-\t c$ is $A$-linear, its right adjoint $f^R=\Map_{C/A}(c,-)$ is lax $A$-linear.
It is exact by stability of $A$ (hence of $C$) and so preserves all colimits if it preserves filtered ones (equivalently, direct sums by stability). 
In order to check the $A$-linearity of $f^R$ it suffices that the map 
$$a \t f^R(c') \r f^R(a \t c')$$
is an isomorphism for all $a \in A$, $c' \in C$.
As both sides preserve colimits in $a$ we may assume $a$ is dualizable by assumption on $A$.
Applying $\Map_{A}(b,-)$, $b\in A$ and using $A$-linearity $f(a^\vee\t b)=a^\vee \t f(b)$ shows that the map is an isomorphism.
So, $f=-\t c$ is internally left adjoint, i.e., $c$ is atomic. 

\textit{Item \eqref{atomic vs compact 2}:}
Let $c\in C$ be compact. 
By assumption, $-\t c\colon A\r C$ preserves compact objects. 
As $A$ is compactly generated, the right adjoint $\Map_{C/A}(c,-)$ preserves filtered colimits.
So, $c$ is atomic by Item \eqref{atomic vs compact 1}. 
\xpf

\subsection{Enriched maps in tensor products}\label{section:infinity-cats-enriched-maps}
Throughout this subsection, let $A \in \CAlg(\PrL)$.

\lemm
\thlabel{tensor adjunctions}
If a functor $f\colon C\to D$ in $\Mod_A(\PrL)$ is internally left adjoint, and $E \in \Mod_A(\PrL)$, then the adjunction $f:C \rightleftarrows  D :f^R$ induces an adjunction:
$$f \t_A \id_E : C \t_A E \rightleftarrows  D \t_A E : f^R \t_A \id_E$$
\xlemm

\pf
First off, $f^R \t_A \id_E$ exists because $f^R$ lies in $\Mod_A(\PrL)$ by assumption.
Applying $- \t_A \id_E$ to the unit and counit map of the adjunction $(f, f^R)$ gives maps $\id \r (f^R \t_A \id_E) \circ (f \t_A \id_E)$ and $(f \t_A \id_E) \circ (f^R \t_A \id_A) \r \id$.
These satisfy the triangle identities, confirming our claim by \cite[Corollary~5.1.14]{Land:Book}.
\xpf

In the following we will use the notation
\eqn
\begin{aligned}
C\x D &\r C\t D \r C\t_A D &&(=\colim_{\Delta^\opp} C\t A^\bullet\t D), \\
(c,d)&\mapsto c\bx d  \mapsto  c\bx_A d 
\end{aligned}
\xeqn  
for objects $c\in C$, $d\in D$.

\prop
\thlabel{proposition:infinity-cats-tensor-product}
Let $C, D \in \Mod_A(\PrL)$.
If $c \in C$ and $d \in D$ are atomic, and $c' \in C, d' \in D$ any objects, then one has
    \begin{equation}
        \label{proposition:infinity-cats-tensor-product:eq1}
        \Map_{C/A}(c, c') \t_A \Map_{D/A}(d, d') \stackrel \cong \lr \Map_{C\t_A D/A}(c\bx_A d, c'\bx_A d').
    \end{equation}
    In particular, $c \bx_A d$ is atomic in $C \t_A D$.

If $C$ and $D$ are atomically generated by $C_0 \subset C^{\at }, D_0 \subset D^{\at }$, then $C \t_A D$ is atomically generated by the image of $C_0 \x D_0$ in $C \t_A D$.
\xprop

\pf
Consider the functors
$$A = A \t_A A \stackrel{- \t c} \lr C = C \t_A A \stackrel{- \boxtimes_A d} \lr C \t_A D.$$
The right adjoint of the composite is $\Map_{C \t_A D/A}(c \bx_A d,-)$.
By \thref{tensor adjunctions}, the right adjoint of the right hand functor is $\id_C \t_A \Map_{D/A}(d,-)$, while the right adjoint of the left hand functor is $\Map_{C/A}(c,-) \t_A \id_A$.
This shows \eqref{proposition:infinity-cats-tensor-product:eq1} and also that $c \bx_A d$ is $A$-atomic.

Further, atomically generated means that the presentable full subcategory of $C$ generated by the objects $a\t c$ for $a\in A$, $c\in C_0$ is $C$ itself, and similarly for $D_0\subset D$.
This implies that the objects $(a \t c) \bx_A (a' \t d)=(a\t a')\t(c\bx_A d)$ for $a,a'\in A$, $c\in C_0$, $d\in D_0$ generate $C \t_A D$ under colimits. 
Since $c\bx_A d$ is atomic, the essential image of $C_0 \x D_0 \r C \t_A D$ atomically generates the target.
\xpf

\subsection{Dualizability of modules}
Throughout, fix $A\in \CAlg(\PrL)$.
The following result characterizes dualizable objects in $\Mod_A(\PrL)$ with respect to $\t_A$ (e.g., this recovers \cite[Proposition~D.7.3.1]{Lurie:SAG} and \cite[Proposition~3.2]{BlumbergGepnerTabuada:UniversalKTheory} if $A=\Spt$):

\theo[{\cite[Observation~1.28, Theorem~1.49]{Ramzi:DualizableCategories}}]\thlabel{dualizability ramzi}
For any object $C\in \Mod_A(\PrL)$, the following are equivalent:
\begin{enumerate}
	\item $C$ is dualizable in $\Mod_A(\PrL)$;
	\item $C$ is a retract of an atomically generated $A$-module.
\end{enumerate}
If, in addition, $C$ is atomically generated, then the enriched Yoneda embedding of \cite{Heine:InftyCategories, Hinich:EnrichedYoneda} induces an equivalence in $\Mod_A(\PrL)$,
\eqn
C \overset{\cong}{\lr} \Fun_A((C^{\textup{at}})^{\opp},A),\;\; c\mapsto \Map_{C/A}(-,c)
\xeqn
with the category of $A$-enriched functors.
\xtheo

\defi\thlabel{definition atomic categories}
The \textit{category of atomically generated $A$-modules} $\Mod_A(\PrL)^{\textup{at}}$ is the non-full subcategory of $\Mod_A(\PrL)$ consisting of atomically generated $A$-modules and $A$-linear functors that preserves atomic objects (equivalently, internal left adjoint functors by \thref{internal left adjoint lemma}).
\xdefi

\exam
If $A=\Spt$, then $\Mod_\Spt(\PrL)^{\textup{at}}=\PrLomegast$ inside $\Mod_\Spt(\PrL)=\PrSt$  (\thref{atomic example}\eqref{atomic example spectra}), i.e, the subcategory of  stable presentable \ii-categories with functors preserving compact objects. 
\xexam

\coro\thlabel{presentable atomics}
The category $\Mod_A(\Pr)^{\textup{at}}$ admits limits and colimits, and the inclusion into $\Mod_A(\PrL)$ preserves colimits. Moreover, the tensor product $\t_A$ preserves this subcategory and turns $\Mod_A(\Pr)^{\textup{at}}$ into a presentably symmetric monoidal category.
\xcoro
\pf
The subcategory is preserved by $\t_A$ (\thref{proposition:infinity-cats-tensor-product}). It is presentable by \cite[Corollary~3.15]{Ramzi:DualizableCategories}, so carries an induced symmetric monoidal structure and admits all (co)limits.
That the tensor product commutes with colimits in both variables follows from the preservation of colimits under the inclusion, which in turn follows from \cite[Theorem~3.7]{Ramzi:DualizableCategories}.
\xpf

\subsection{Enforcing conservativity}
\label{sect--enforcing conservativity}
The results of this section are used in \S\ref{section:cohomological-motives} to construct cohomological motives.
Throughout, we denote by $\PrSt=\Mod_{\Spt}(\PrL)$ the category of stable presentable categories (and continuous functors).

Recall from \cite[Section 5]{BlumbergGepnerTabuada:UniversalKTheory} that a sequence $$C\to D\to E$$ in $\PrSt$ is called \textit{exact} if it is a fiber and cofiber sequence in $\PrSt$.
Equivalently, $D\to E$ is a localization at all arrows in $E$ whose (co)fiber lies in $D$ with kernel $C$. 
Given an exact sequence $C\to D\to E$ in $\PrSt$, we call $E=D/C$ the \textit{Verdier quotient of $D$ by $C$}.

For a map $f\colon C\r D$ in $\PrSt$ we denote by $\ker f$ the full subcategory of objects $c\in C$ with $f(c)\simeq 0$.
Passing to the Verdier quotient induces a diagram
\eqn\label{conservativity factorization}
\xymatrix{C \ar[rr]^f \ar[dr]_{f'} & & D \\ & C / \ker f \ar[ur]_{f''}}
\xeqn
in which $f'$ is a localization and $f''$ is conservative (by stability it is enough
to check that an object $c \in C$ is zero in $C / \ker f$ if and only if $f''(c) = 0$, but this holds precisely by construction).

\lemm
\thlabel{factorization}
There is a lax symmetric monoidal functor
$$\Fun(\Delta^1, \PrSt) \r \Fun(\Delta^2, \PrSt)$$
sending a map $f\colon C \r D$ in $\PrSt$ to the diagram \eqref{conservativity factorization}.
\xlemm
\pf
Let $I$ denote the (ordinary) category $\{0 \r 1\}$.
The functor $\Fun(I, \PrSt) \r \Fun(\{0 \r 1 \gets 0' \}, \PrSt)$ sending some functor $C \r D$ to the diagram $C \r D \gets 0$ is symmetric monoidal.
Next, the functor $\Fun(\{0 \r 1 \gets 0' \}, \PrSt) \r \Fun(I^2, \PrSt)$ sending a diagram $C \r D \gets D'$ to the pullback square involving $C \x_D D'$ is lax symmetric monoidal.
The functor sending that pullback diagram to $C \x_D D' \r C$ is again symmetric monoidal.

Consider the symmetric monoidal structure on $I$ given by the minimum function, i.e., $0 \t 0 = 0 \t 1 = 1 \t 0 = 0$, $1 \t 1 = 1$.
We denote the resulting Day convolution product on the functor category $\Fun(I, \PrSt)$ by $\diamond$. 
Concretely, given two arrows $r\colon C \r D$, $r'\colon C' \r D'$, we have $r \diamond r' \colon D \t C' \sqcup_{C \t C'} C \t D' \r D \t D'$.
Recall that $-\diamond-$ is the composite
$$-\diamond- \colon \Fun(I, \PrSt)^2 \stackrel{\t} \lr \Fun(I^2, \PrSt) \stackrel{\min^*} \lr \Fun(I, \PrSt),$$
where $\min^*$ is the left adjoint of the restriction functor $\min_*$.
This induces natural maps
$$(r \t r') \diamond (s \t s') \r (r \diamond s) \t (r' \diamond s').$$
These assemble to the structure of a oplax symmetric monoidal functor for $\diamond$.
In particular $r \diamond (\Spt \r 0) = (\cofib r \r 0)$, so composing all the above-mentioned lax symmetric monoidal functors gives rise to a lax symmetric monoidal functor $\Fun(I, \PrSt) \r \PrSt$, $(C\stackrel f \r D) \mapsto C / \ker f$.

There is another symmetric monoidal functor $$\Fun(\Delta^1, \PrSt) \r \Fun(\Delta^2, \Fun(\Delta^1, \PrSt)),$$ sending an arrow $C \r D$ in $\PrSt$ to $\id_C \r (C \r D) \r \id_D$.
Composing all these gives the requested functor, mapping a map of commutative algebra objects in $\PrSt$ to 
$C \stackrel{f'} \r C / \ker f \stackrel {f''} \r D$.
\xpf

Let $S$ be a simplicial set. 
Similarly to \cite[Definition~4.7.4.16]{Lurie:HA}, we denote by
\begin{equation}
\Fun'(S, \Fun(\Delta^n, \PrLSt)) \subset \Fun(S, \Fun(\Delta^n, \PrSt)) = \Fun(S \x \Delta^n, \PrLSt)
\label{Fun'}
\end{equation}
the full subcategory spanned by functors $f \colon S \x \Delta^n \r \PrLSt$ such that for each map $\sigma \colon s \r \ol s$ in $S$, and each map $\tau \colon [a] \r [\ol a]$ in $\Delta^n$, the diagram
\begin{equation}
    \xymatrix{
f(s, [a]) \ar[r] \ar[d] & f(s, [\ol a]) \ar[d] \\
f(\ol s, [a]) \ar[r] & f(\ol s, [\ol a])
}
\label{adjointable}
\end{equation}
is \emph{strongly} right adjointable, by which we mean that the vertical functors admit continuous right adjoints, and these right adjoints commute with the displayed horizontal maps.

\lemm
\thlabel{factorization and right adjoints}
The factorization functor supplied by \thref{factorization} restricts to a functor
$$\Fun'(S, \Fun(\Delta^1, \PrLSt)) \r \Fun'(S, \Fun(\Delta^2, \PrLSt)).$$
\xlemm

\pf
Let $\varphi \in \Fun'(S, \Fun(\Delta^1, \PrLSt))$ and let $f \colon S \r \Fun(\Delta^2, \PrLSt)$ be its factorization. 
We write $f_k(s) = f(s, k)$, so that $f_0(s) \r f_2(s)$ is the same as $\varphi(s)$.
For a map $\sigma \colon s \r \ol s$ in $S$, we have to check that the following diagram is (vertically) strongly right adjointable.
\begin{equation}
    \xymatrix{
f_0(s) \ar[r]^(.3){a'} \ar[d]_{\sigma_0} & 
f_1(s) = f_0(s) / \ker a \ar[r]^(.7){a''} \ar[d]_{\sigma_1} & 
f_2(s) \ar[d]_{\sigma_2} \\
f_0(\ol s) \ar[r]^(.3){\ol a'} \ar@/_/@{.>}[u]_{\sigma_0^R} & 
f_1(\ol s) = f_0(\ol s) / \ker \ol a \ar[r]^(.7){\ol a''} \ar@/_/@{.>}[u]_{v} & 
f_2(\ol s). \ar@/_/@{.>}[u]_{\sigma_2^R} 
}
\label{vertically right adjointable}
\end{equation}
By assumption, $\sigma_0^R$ is continuous, and its restriction to $\ker \ol a$ factors over $\ker a$ (since the outer diagram is adjointable). 
By the universal property of the Verdier quotient in $\PrLSt$, there is a continuous functor
$$v \colon f_0(\ol s) / \ker \ol a \r f_0(s) / \ker a$$
as displayed, making the left square (with the two curved functors) commute. By construction, there are natural maps $\id \r \sigma_1 \circ v$, $\id \r v \circ \sigma_1$, and they satisfy the triangle identities.
Thus, $v$ is a right adjoint of $\sigma_1$.
This confirms the adjointability of the left hand square.
Using that $\ol a'$ is a localization (marked * below), the adjointability of the right square then follows from the one of the total diagram:
$$
a'' \sigma_1^R \stackrel * = a'' \sigma_1^R \ol a' \ol a'^R =
a'' a' \sigma_0^R \ol a'^R =
a \sigma_0^R \ol a'^R = \sigma_2^R \ol a \ol a'^R
= \sigma_2^R \ol a'' \ol a' \ol a'^R \stackrel * = \sigma_2^R \ol a''.
$$
\xpf

\rema
In the proof, the existence of a right adjoint of $\sigma_1$ is a priori clear, since $\sigma_1$ is in $\PrLSt$. The benefit of assuming that the outer diagram is \emph{strongly} vertically adjointable is this explicit description of the right adjoint $\sigma_1^R$.
\xrema

\section{An abstract K\"unneth formalism}\label{section:abstract-kuenneth}
\label{section:abstract-kuenneth-set-up}
In this section, we discuss the categorical K\"unneth formula on an abstract categorical level. 
More precisely, for a category $C$ of geometric objects (e.g., finite type schemes over a field) and an abstract $6$-functor formalism $\D$ on $C$ (e.g., a motivic sheaf formalism) we give criteria for which the exterior product functor \eqref{functor bx} is fully faithful (\thref{proposition:abstract-fully-faithful}), respectively an equivalence (\thref{proposition:essentially-surjective,Künneth inddualizable}).

\subsection{Setting} We fix an \ii-category $C$ that has finite limits, including a final object $\star \in C$.
Then, $C^\opp$ admits finite coproducts, hence carries the coCartesian symmetric monoidal structure \cite[2.4.3.3]{Lurie:HA}.
We denote by $\Corr(C)$ the symmetric monoidal \ii-category of correspondences from \cite[Definition~A.5.4]{Mann:SixFunctors} (cf.~\cite[Section~5.1]{Chowdhury:SixIII}) with \textit{all} morphisms in $C$. 
For the purposes of this paper, though, it is also possible to use the definition of \cite{GaitsgoryRozenblyum:StudyI} or the construction in \cite[Section~4.6]{AyoubGallauerVezzani:RigidAnalyticMotives}.

We fix a lax symmetric monoidal functor
\eqn\label{equation:definition-kuenneth-setting-naive}
\D  \colon \Corr(C) \to \PrL.
\xeqn
Any $X \in C$ is naturally a (commutative) algebra object over the final object $\star$. 
So, $\D$ induces a lax symmetric monoidal functor
 \begin{equation}
    \label{equation:definition-kuenneth-setting}
    \D\colon \Corr(C) =\Mod_\star(\Corr(C)) \r \Mod_{\D(\star)}(\PrL). \end{equation}
In particular, for objects $X, Y \in C$, the datum of $\D$ includes a functor
\begin{equation}
    -\bx- := - \bx_{X, Y} - : \D(X) \t_{\D(\star)} \D(Y) \r \D(X \x Y).
    \label{functor bx}
\end{equation}
Because of the examples in the following sections, we refer to this functor as the \emph{exterior product}.

\defi\thlabel{definition-categorical-kuenneth}
The \textit{categorical Künneth formula holds (for $\D$ and objects $X, Y \in C$)} if \eqref{functor bx} is an equivalence.
\xdefi

Throughout, the structural map to the final object is denoted $p_X\colon X\r \star$.

\nota
\thlabel{notation functors}
We introduce the following notation pertaining to an abstract $6$-functor formalism:
\begin{enumerate}
    \item \label{notation:kuenneth-setting:it2}
        By \cite[Theorem~2.4.3.18]{Lurie:HA}, any $X \in C$ is naturally a commutative coalgebra in $C$ (by means of the diagonal $\Delta_X : X \r X \x X$). 
        Thus $\D(X)$ is naturally a commutative algebra object in $\PrL$, whose monoidal structure we denote by $\t_X$ or just $\t$.
        Concretely, for $c, d \in \D(X)$, $c \t d = \Delta_X^*(c \boxtimes d)$.
        Being presentable, $\D(X)$ is then a closed symmetric monoidal \ii-category, whose inner hom we denote by $\IHom_X(-,-)$.
               \item
       \label{f* f! notation}
    For a map $f : X \r Y$ in $C$, we write 
        \eqnarra
            f^*:=\D(X \stackrel \id \gets X \stackrel f \r Y)\colon \D(Y)\to \D(X), \label{equation:notation-pullback} \cr
            f_! := \D(Y \stackrel f \gets X \stackrel \id \r X) : \D(X) \r \D(Y). \nonumber
        \xeqnarra
    Both functors are $\D(\star)$-linear (in fact, $\D(Y)$-linear).
    Being functors in $\PrL$, they admit right adjoints $f_* := (f^*)^R$ and $f^! := (f_!)^R$, which are lax $\D(\star)$-linear.
    In these terms, the functor \eqref{equation:definition-kuenneth-setting-naive} encodes that the formation of $f_!$, $f^*$, $f^!$, $f_*$ is functorial in $f$ and satisfies base change and the projection formula \cite[Proposition~A.5.8]{Mann:SixFunctors}.
     \item
    \label{only 4 functors}
    Restricting the functor $\D$ along $C^\opp \r \Corr(C)$ gives a lax symmetric monoidal functor $C^\opp \r \PrL$, or equivalently by \cite[Theorem~2.4.3.18]{Lurie:HA} a functor
    $$C^\opp \r \CAlg(\PrL), X\mapsto \D(X), f\mapsto f^*.$$
    A number of abstract results below only use this part of the functoriality.
\end{enumerate}
\xnota

The following lemma relates the inner homs with the enriched maps from \eqref{enriched mapping functor}:

\lemm\thlabel{lemma:inner-hom}
For each map $f\colon X\to Y$ in $C$, one has
\eqn
	\Map_{\D(X)/\D(Y)}(-,-)=f_*\IHom_{X}(-,-)
\xeqn
as functors $\D(X)^\opp\x \D(X)\r \D(Y)$.
\xlemm
\pf
The left hand side is well-defined by \thref{notation functors}\eqref{only 4 functors}.
One checks that $f_*\IHom_X(-,-)$ satisfies the defining property of the enriched mapping functor.
\xpf

Recall from \thref{presentable atomics} that the subcategory $\Mod_{\D(\star)}(\PrL)^{\textup{at}}$ of atomically generated $\D(\star)$-modules in $\Mod_{\D(\star)}(\PrL)$ is stable under $\t_{\D(\star)}$ and carries the induced symmetric monoidal structure.
We often impose the following assumption which is satisfied in all our examples.

\assu\thlabel{useful assumption} 
The functor $\D\colon \Corr(C) \r \Mod_{\D(\star)}(\PrL)$ takes values (as a lax symmetric monoidal functor) in atomically generated $\D(\star)$-modules $\Mod_{\D(\star)}(\PrL)^{\textup{at}}$. 
More concretely, this means that 
\begin{enumerate}
	\item \label{useful assumption objects} 
		$\D(X)$ is atomically generated for all objects $X\in C$, 
	\item  \label{useful assumption maps}
		$f^*, f_!$ preserve atomic objects for all maps $f\colon X\r Y$ in $C$, and
	\item \label{useful assumption hom}
		the exterior product \eqref{functor bx} preserves atomic objects for all $X,Y\in C$.
\end{enumerate}
If $\D$ is only defined on $C^\opp$ (cf.~\thref{notation functors}\eqref{only 4 functors}), then \eqref{useful assumption maps} is imposed only for $f^*$.
\xassu

\rema
Under \thref{useful assumption} the functor $-\t_XM\colon \D(X)\r \D(X)$ is an internal left adjoint for any $M\in \D(X)^{\textup{at}/\D(\star)}$, i.e., $\IHom_X(M,-)$ is continuous and $\D(\star)$-linear.
\xrema

\defi\thlabel{definition indbl}
Let $\D(X)^\indbl\subset \D(X)$ be the subcategory of \textit{ind-dualizable objects}, i.e., the presentable full subcategory generated by dualizable objects. 
This induces a lax symmetric monoidal subfunctor $\D(-)^\indbl\colon C^\opp\r \PrL$ of the restriction $\D|_{C^\opp}$ (cf.~\thref{notation functors}\eqref{only 4 functors}).
\xdefi

The next lemma shows that under a mild additional assumption $\D(-)^\indbl$ takes values in atomically generated modules if $\D$ does:

\lemm\thlabel{enriched maps ind dualizable}
Suppose $\D$ satisfies \thref{useful assumption}.
In addition, suppose $\D(\star)$ is generated under colimits by dualizable objects, i.e., $\D(\star)^\indbl=\D(\star)$.
Then, for all $X\in C$, the dualizable objects $M\in \D(X)$ are $\D(\star)$-atomic and
\[
\Map_{\D(X)^\indbl/\D(\star)}(M,-)=p_{X,*}\IHom_{X}(M,-)
\]
as functors $\D(X)^\indbl\r \D(\star)^\indbl$.
In particular, $\D(-)^\indbl$ satisfies \thref{useful assumption}.
\xlemm
\pf
The inclusion $\D(X)^\indbl\subset \D(X)$ is symmetric monoidal. 
Since the functor $\IHom_X(M,-)=M^\vee\t_X -$ commutes with colimits and preserves dualizable objects its restriction to $\D(X)^\indbl$ computes the inner hom in this subcategory. 
Using the assumption $\D(\star)^\indbl=\D(\star)$, we see that the adjunction $p_X^*:\D(\star)\leftrightarrows \D(X): p_{X,*}$ restricts to the subcategories of ind-dualizable objects. 
Since $p_{X,*}$ is continuous and $\D(\star)$-linear under \thref{useful assumption}, the lemma follows from \thref{lemma:inner-hom} applied to $\D(-)^\indbl$.
\xpf

\subsection{Full faithfulness}
The next result shows that the full faithfulness of \eqref{functor bx} should be regarded as a replacement of the Künneth formula for exterior products of objects:

\prop\thlabel{proposition:abstract-fully-faithful}
Fix $X,Y \in C$. 
Suppose $\D(X), \D(Y)$ are $\D(\star)$-atomically generated by subcategories $\D_0(X), \D_0(Y)$ respectively.  
Then, the exterior product \eqref{functor bx} is fully faithful if and only if the map
\eqn\label{proposition:fully-faithful:eq1}
p_{X,*}{\IHom}_X(M,M')\t p_{Y,*}{\IHom}_Y(N,N'){\r} p_{X\times Y,*}\IHom_{X\x Y}(M\bx N,M'\bx N')
\xeqn
is an isomorphism in $\D(\star)$ for all $M, M'\in \D_0(X)$, $N,N'\in \D_0(Y)$.
\xprop
\pf
By \thref{proposition:infinity-cats-tensor-product}, $\D(X) \t_{\D(\star)} \D(Y)$ is atomically generated by the atomic objects $M \boxtimes N$.
Our claim now follows the computation of enriched maps in \thref{lemma:inner-hom} and \thref{internal left adjoint lemma}.
\xpf

\coro\thlabel{proposition:abstract-fully-faithful indbl}
Suppose $\D$ satisfies \thref{useful assumption} and $\D(\star)^\indbl=\D(\star)$.
Then, the exterior product for $\D(-)^\indbl$ and $X,Y$ is fully faithful if and only if \eqref{proposition:fully-faithful:eq1} is an isomorphism for all dualizable objects $M, M'\in \D(X)$, $N,N'\in \D(Y)$.
\xcoro
\pf
This follows from \thref{enriched maps ind dualizable}.
\xpf

\subsection{Essential surjectivity}\label{section:essentially-surjective}

\subsubsection{Conservativity of $\bx^R$ in a descent-type situation}\label{section:conservativity-six}

Given that a fully faithful functor $f$ is an equivalence if and only if it admits a \emph{conservative} right adjoint $f^R$, we provide a condition that ensures the conservativity of the functor
$$\boxtimes^R : \D(X \x Y) \r \D(X) \t_{\D(\star)} \D(Y),$$
right adjoint to $\bx$ from \eqref{functor bx}.
Note that $\bx^R$ exists since $\bx$ lies in $\PrL$.
The following lemma is applied in the proof of \thref{Künneth DMf indlis} to perform a descent argument.

\prop
\thlabel{pushforward conservative}
\thlabel{Kuenneth constructible}
Suppose that $\D$ satisfies \thref{useful assumption} and that the exterior product functor \eqref{functor bx} is fully faithful for all pairs of objects in $C$.
In addition, suppose that there is a family of correspondences
$$\{Y \stackrel {g_i} \gets Z_i \stackrel{f_i} \r Y_i\}_{i \in I}$$
in $C$ such that 
$$\prod_{i\in I}(\id_X \x f_i)_* (\id_X \x g_i)^! : \D(X \x Y) \r \prod_{i\in I} \D(X \x Y_i)$$
is conservative.
Then, the categorical Künneth formula for all pairs $(X, Y_i)$, $i\in I$ implies the one for $(X, Y)$.
\xprop

\pf
By the lax monoidality of $\D$, the following diagram commutes
$$\xymatrix{
\D(X) \t \D(Y_i) \ar[d]_{\id \t g_{i!} f_i^*} \ar[rr]^{\bx_{X, Y_i}} & &
\D(X \x Y_i) \ar[d]_{(\id \x g_i)_! (\id \x f_i)^* } \\
\D(X) \t \D(Y) \ar[rr]^{\bx_{X,Y}} & &
\D(X \x Y),
}$$
where $\t:=\t_{\D(\star)}$ and $\bx_{-,-}$ denotes the exterior product functor \eqref{functor bx} for the respective objects. 
Passing to right adjoints (where $\id \t g_i^! f_{i*}$ is the right adjoint of $\id \t g_{i!} f_i^*$ by \thref{tensor adjunctions} using \thref{useful assumption}) gives a commutative diagram and so
$$\prod_{i\in I} \bx_{X, Y_i}^R \circ (\id \x g_i)^! (\id \x f_i)_* = \prod_{i\in I} (\id \t g_i^! f_{i,*}) \circ \bx_{X, Y}^R.$$
By assumption $\bx_{X, Y_i}^R$ is an equivalence and $\prod_{i\in I} (\id \x g_i)^! (\id \x f_i)_*$ is conservative. 
Hence, $\bx_{X,Y}^R$ is conservative as well.
\xpf

\subsubsection{An explicit construction of a preimage}
The following construction, which is analogous to Gaitsgory--Kazhdan--Rozenblyum--Varshavsky \cite[A.5.7]{GaitsgoryEtAl:Toy}, constructs a candidate for a preimage of a given object in $\D(X \x Y)$. 
It works for functors $\D$ just defined on $C^\opp$ (cf.~\thref{notation functors}\refit{only 4 functors}).

\const\thlabel{construction:essentially-surjective}
Let $X, Y\in C$. 
Assume that $p_X\colon X\r \star$ has a section $x$ and consider the following Cartesian diagram:
$$\xymatrix{
Y \ar[r]^(.4){x \x \id} \ar[d]^{p_Y} & 
X \x Y \ar[d]^{\pr_X} \ar[r]^(.6){\pr_Y} & 
Y \ar[d]^{p_{Y}} \\ 
\star \ar[r]^x & 
X \ar[r]^{p_X} & 
\star}$$
For an object $P\in \D(X\x Y)$, define the following objects:
\eqn
\begin{aligned}
& N :=(x\times\id_Y)^*P & \in  \D(Y) \\
& M :=\pr_{X,*}\underline{\Hom}_{X\x Y}(\pr_{Y}^*N, P)&\in \D(X)\\
& E :=\pr_X^*\pr_{X,*}\underline{\End}_{X\x Y}(\pr_{Y}^*N)&\in \D(X\x Y)
\end{aligned}
\xeqn
\xconst

\lemm\thlabel{lemma:essentially-surjective}
With the notation of \thref{construction:essentially-surjective}, the following holds:
\begin{enumerate}
	\item\label{lemma:essentially-surjective:it1}	
	The object $E$ has the structure of an $\bbE_1$-algebra in $\D(X\x Y)$. In addition, $\pr_X^*M$ is a right $E$-module, while $pr_Y^* N$ is a left $E$-module.
	\item\label{lemma:essentially-surjective:it2}
	The evaluation map $M\bx N=\pr_X^*M\t\pr_{Y}^*N\to P$, $f\otimes n\mapsto f(n)$ has an $E$-linear structure and hence factors through a map
	\eqn\label{lemma:essentially-surjective:eq1}	
	M\bx_EN:=\colim\pr_X^*M\t E^\bullet \t\pr_{Y}^*N \to P,
	\xeqn
	where $\t:=\t_{X\x Y}$ and the colimit is the two-sided bar construction from \cite[Construction~4.4.2.7]{Lurie:HA}.
	\item\label{lemma:essentially-surjective:it3}
    Suppose that the exterior product functor \eqref{functor bx} is fully faithful. Then $M \bx_E N$ lies in its essential image.
\end{enumerate}
\xlemm
\pf
\eqref{lemma:essentially-surjective:it1}: 
Like all $*$-pullback functors, $\pr_X^*$ is symmetric monoidal (\thref{notation functors}\eqref{only 4 functors}), so its right adjoint $\pr_{X*}$ is lax symmetric monoidal.
Since $\underline{\End}_{X\x Y}(\pr_{Y}^*M_{Y})$ is an algebra object in $\D(X\x Y)$ with a left action on $\pr_{Y}^*N$ \cite[Section~4.7.1]{Lurie:HA}, the object $E$ comes with the desired structures by precomposing with the counit $\pr_X^*\pr_{X,*}\to \id_{\D(X\x Y)}$.

\eqref{lemma:essentially-surjective:it2}:
The $E$-linearity is immediate.

\eqref{lemma:essentially-surjective:it3}:
By definition, the elementary tensors $M\bx N$ and $E=\pr_{X,*}\underline{\End}(\pr_{Y}^*N)\bx\unit_{Y}$ are in the essential image of $\D(X) \x \D(Y) \r \D(X) \t_{\D(\star)} \D(Y)$. 
Hence, so are $M\bx_EN$ because the exterior product functor \eqref{functor bx} is fully faithful, symmetric monoidal and preserves colimits.
\xpf

The next statement gives a condition that ensures that the candidate \eqref{lemma:essentially-surjective:eq1} is indeed a preimage.

\prop\thlabel{proposition:essentially-surjective}
Let $\D' \colon C^\opp \r \PrL$ be a lax symmetric monoidal functor that satisfies \thref{useful assumption} and whose exterior product functor \eqref{functor bx} is fully faithful for all objects in $C$.
Fix $X,Y\in C$ such that $p_X\colon X\r \star$ has a section $x$.
Suppose that for all $\D'(\star)$-atomic objects $P \in \D'(X \x Y)$ the following conditions are satisfied (cf.~\thref{construction:essentially-surjective} for notation):
\begin{enumerate}
	\item \label{proposition:essentially-surjective:it1}
		The evaluation of the following natural transformation (of functors $\D'(X\x Y)\to \D'(\star)$),
		\eqn  \label{proposition:essentially-surjective:eq1}
		x^*\pr_{X,*}\r p_{Y,*}(x\x \id_Y)^* 
		\xeqn
		at the objects $\IHom_{X\x Y}(\pr_Y^*N,P)$, $\underline{\End}_{X\x Y}(\pr_{Y}^*N)$ is an isomorphism.
	\item \label{proposition:essentially-surjective:it2}
		The evaluation of the following natural transformation (of functors $\D'(X\x Y)^\opp\x \D'(X\x Y)\to \D'(Y)$),
		\eqn  \label{proposition:essentially-surjective:eq2}
		(x\x \id_Y)^*\IHom_{X\x Y}(-,-)\to \IHom_{Y}(-,-)\circ (x\x \id_Y)^*
		\xeqn
		at the objects $(\pr_Y^*N,P)$ and $(\pr_Y^*N,\pr_Y^*N)$ is an isomorphism.
    \item \label{pullback conservative}The functor
\eqn \label{proposition:essentially-surjective:eq3}
 (x\x \id_Y)^*\colon \D'(X\x Y)\to \D'(Y)
\xeqn
is conservative.
\end{enumerate}
Then, the categorical K\"unneth formula holds for the objects $X, Y$ (\thref{definition-categorical-kuenneth}).
\xprop
\pf
First off, $\D'(X\x Y)$ is $\D'(\star)$-atomically generated by \thref{useful assumption}.
So, using Item \eqref{pullback conservative} it is enough to show that the map $M\bx_E N\to P$ defined in \eqref{lemma:essentially-surjective:eq1} is an equivalence after applying $(x\times\id_Y)^*$.
We compute:
\begin{equation}\label{equiblah}
(x\times\id_Y)^*\pr_X^*M\otimes_{(x\times\id_Y)^*E}N=(x\times\id_Y)^*(M\bx_EN) \to (x\times \id_Y)^*P=N
\end{equation}
Using Items~\eqref{proposition:essentially-surjective:it1} and \eqref{proposition:essentially-surjective:it2} gives 
$$
\begin{aligned}
& (x\times\id_Y)^*\pr_X^*M=p_Y^*p_{Y,*}(x\times \id_Y)^*\underline{\Hom}(\pr_Y^*N,P)=p_Y^*p_{Y,*}\underline{\End}(N) \\
& (x\times\id_Y)^*E=p_Y^*p_{Y,*}(x\times \id_Y)^*\underline{\End}(\pr_Y^*N)=p_Y^*p_{Y,*}\underline{\End}(N)
\end{aligned}
$$
compatibly with the right $E$-module structure on $\pr_X^*M$ and where the inner homs are computed in $\D'(X\x Y)$.
So, the map \eqref{equiblah} becomes the identity map 
\eqn
N=p_Y^*p_{Y,*}\underline{\End}(N)\otimes_{p_Y^*p_{Y,*}\underline{\End}(N)}N \to N,
\xeqn
which implies the proposition.
\xpf

In the context of motivic sheaves, the most critical assumption is the conservativity of $(x\x \id_Y)^*$; it provides a motivation for the study of cohomological motives in Section \ref{section:cohomological-motives} below.
The other two assumptions belong to the more classical realm of six functors, as the following statement shows (cf.~\thref{* pushforward remark,purity remark}).

\prop
\thlabel{Künneth inddualizable}
Let $\D\colon \Corr(C) \r \PrL$ be a lax symmetric monoidal functor. 
Then, Items
\refit{proposition:essentially-surjective:it1} and \refit{proposition:essentially-surjective:it2} in \thref{proposition:essentially-surjective} are satisfied for $\D'=\D(-)^\indbl$ provided that 
\begin{enumerate}
    \item \label{pr* dualizable} $\pr_{X,*} \colon \D(X \x Y) \r \D(X)$ preserves dualizable objects (i.e., restricts to $\D(-)^\indbl$),
    \item \label{x pure} the object $x^! 1_X$ is $\t$-invertible, and
    \item \label{*! exchange} the natural map 
    \begin{equation}
     p_Y^* x^! 1_X \r (x \x \id_Y)^! \pr_X^* 1_X \label{BC * ! unit}
    \end{equation} 
        is an isomorphism.
        \end{enumerate}
\xprop
\pf
\textit{Item~\refit{proposition:essentially-surjective:it2} in \thref{proposition:essentially-surjective}:} The map in
\eqref{proposition:essentially-surjective:eq2} is an isomorphism if $P$ and therefore $\pr_Y^* N$ are dualizable in $\D(X \x Y)$, since then
\eqn
\IHom_{X\x Y}(\pr_Y^*N,-)=(\pr_Y^*N)^\vee\t_{X\x Y}(-).
\xeqn

\textit{Item~\refit{proposition:essentially-surjective:it1} in \thref{proposition:essentially-surjective}:} 
We need to show that the natural map
$$\text{BC}^*_{*} \colon x^* \pr_{X,*} \r p_{Y,*} (x \x \id_Y)^*$$ 
is an isomorphism when evaluated at a dualizable object $F\in \D(X \x Y)$.

We use the following facts which are readily checked. 
For any map $s\colon T \r S$ in $C$, the natural map (adjoint to the projection formula)
$$\delta\colon s^* A \t_T s^! 1 \r s^! A.$$ 
is an isomorphism for any dualizable object $A\in \D(S)$.
Further, the monoidality of $s^*$ yields (by passing to adjoints) a natural map
$$\alpha\colon s_* B \t_S A \r s_*(B \t_T s^* A)$$
for any $B\in \D(T)$.
Again, for $A\in \D(A)$ dualizable, this map is an isomorphism.

With these preparations, we consider the following commutative diagram \cite[Lemma~1.1.2]{Saito:Cotangent} (see also Proposition 1.1.8, Part (2) there), where the top left horizontal map is $\text{BC}^*_* \t \id_{x^! 1}$, while the bottom left map is the base change isomorphism $\text{BC}_*^!$ that is part of the correspondence functoriality of $\D$:
$$\xymatrix{
x^* \pr_{X*} F \t x^! 1 \ar[r]^(.45){\text{BC}^*_* \t \id} \ar[d]^{\cong \text{ by \eqref{pr* dualizable}}}_{\delta} & 
p_{Y*} (\id \x x)^* F \t x^! 1 \ar[r]^{\cong \text{ by \eqref{x pure}}}_\alpha &
p_{Y*} ((\id \x x)^* F \t p_Y^* x^! 1) \ar[d]^{\cong \text{ by \eqref{*! exchange}}}_{\eqref{BC * ! unit}} \\
x^! \pr_{X*} F \ar[r]^\cong_{\text{BC}^!_*} &
p_{Y*} (\id \x x)^! F &
p_{Y*} ((\id \x x)^* F \t (\id \x x)^! 1) \ar[l]^(.6)\cong_(.6){\delta}
}$$
Thus, under our conditions \eqref{pr* dualizable}--\eqref{*! exchange}, the top left map is an isomorphism. 
So, tensoring with the $\t$-inverse of $x^! 1$, we obtain that $\text{BC}^*_*$ evaluated at $F$ is an isomorphism, as requested.
\xpf

\rema\thlabel{* pushforward remark}
The question when $*$-pushforward preserves dualizable objects, or commutes with $*$-pullbacks has been considered in various places including 
\cite[Th.~finitude, Appendix, 5.1.3]{SGA412} (in the situation of tamely ramified étale sheaves), 
\cite[Corollaire~2.1.2]{Laumon:Semi-Continuite} (for relative curves, in terms of the Swan conductor) and
\cite[Proposition~2.4.4]{Cisinski:Cohomological} (an abstract criterion for h-motives).
\xrema

\rema\thlabel{purity remark}
If $\D$ is one of the $6$-functor formalisms as in our examples below, then $x\colon \star \r X$ will be the inclusion of a point into the regular locus of the variety $X$.
In this case, Item \eqref{x pure} and \eqref{*! exchange} of \thref{Künneth inddualizable} hold by purity and the compatibility of the fundamental class with respect to transversal pullback \cite[Corollary~1.3.22, Proposition~1.3.21(b)]{BondarkoDeglise:Homotopy-t-structure}.
\xrema

\section{A conjecture for étale motives}\label{section:algebraic-motives}

Let $\Lambda\in \CAlg(\Spt)$ be a connective, commutative ring spectrum, and abbreviate $\Mod_\Lambda=\Mod_\Lambda(\Spt)$.
We fix a field $k$ and work with the category of finite type $k$-schemes $\Sch^\ft_k$.
For such a scheme $X$, let $\Sm_{/X}$ be the full subcategory of $\Sch_{/X}$, the category of schemes over $X$, consisting of $X$-smooth maps.
Throughout, we work with the following model of the category of \textit{étale motives} \cite{Ayoub:Realisation, AyoubGallauerVezzani:RigidAnalyticMotives} defined by the adjunctions between $\Lambda$-linear, presentably symmetric monoidal \ii-categories:  
\eqn\label{equation:SH-adjunctions}
	{\HypShv}_\et({\Sm}_{/X},\Lambda) \overset{{L_{\mathbb A^1}}}{\underset{\textrm{incl}} \rightleftarrows} 
	\DA_{\et,\eff}(X,\Lambda) \overset{{\Sigma_{\Gm}}}{\underset{\Omega_{\Gm}} \rightleftarrows}
	\DA_\et(X,\Lambda) 	
\xeqn
Here, $\HypShv_\et(\Sm_{/X},\Lambda)$ is the \ii-category of étale hypersheaves on $\Sm_{/X}$ with values in $\Mod_\Lambda$ and $\DA_{\et,\eff}(X,\Lambda)$ its full subcategory of $\A^1$-invariant sheaves. 
The functor $\Sigma_{\bbG_m}$ is the initial functor in $\CAlg(\PrL)$ sending $1_X(1):=\Gm[-1]$ to an invertible object.

\exam\thlabel{example:coefficients}
The following special cases are of interest to us:
\begin{enumerate}
    \item \label{example:coefficients:it3}
    If $\Lambda$ is discrete (i.e., $\Lambda=\pi_0\Lambda$), then \eqref{equation:definition-SH} is Ayoub's category of étale motives \cite{Ayoub:Realisation}. In particular, for $\Lambda=\Q$ one gets the category of Beilinson motives \cite[Theorem~5.2.2, Corollary~5.5.7]{CisinskiDeglise:Etale}.
    \item 
    \label{Z/ell}
    If $\Lambda = \Z/n$ with $n\in \bbZ$ invertible in $k$,
	then rigidity \cite[Théorème~4.1]{Ayoub:Realisation} (see also \cite[Theorem~6.6]{Bachmann:Rigidity} and \cite[Theorem~3.1]{BachmannHoyois:Rigidity}) implies an equivalence
$$\HypShv(X_\et,\Z/n)\cong \DA_{\et}(X,\Z/n),$$
where at the left we have the \ii-category of étale hypersheaves of $\Z/n$-modules.
In addition, note that $\DA_{\et}(X,\Z/p)=0$ if $\Char(k)=p>0$ by \cite[Corollary~4.5.3]{CisinskiDeglise:Etale}, \cite[Theorem~A.1]{BachmannHoyois:Rigidity}.
Also, recall from 
\cite[Theorem~2.1.2.2, Remark~2.1.2.1]{Lurie:SAG} that its homotopy category is the classical unbounded derived category of étale $\Z/n$-sheaves.
    \item \label{example:coefficients:it1} 
            We may replace the étale topology above by the Nisnevich topology. In this case, if we choose $\Lambda=\bbS$ (the sphere spectrum), then \eqref{equation:SH-adjunctions} defines the stable $\A^1$-homotopy category $\SH$.
\end{enumerate}
\xexam

The construction upgrades to a lax symmetric monoidal functor 
\eqn\label{equation:definition-SH}
\DA := \DA_\et(-,\Lambda)\colon \Corr(\Sch^\ft_k) \to \PrL
\xeqn
where $\Corr(\Sch^\ft_k)$ denotes the category of correspondences:
Briefly, one starts with the functor $(\Sch^{\ft}_k,\textup{sep})^\opp\r \CAlg(\PrL), X\mapsto \DA(X), f\mapsto f^*$ constructed in \cite[\S9.1]{Robalo:Theorie} and restricted to separated morphisms (hence, every morphisms can be factored as an open immersion followed by a proper morphism). 
Then, one checks the conditions in \cite[Proposition~A.5.10]{Mann:SixFunctors} using \cite[\S6]{Ayoub:Realisation}, so this functor extends to the lax symmetric monoidal functor $ \Corr(\Sch^{\ft}_k,\textup{sep}) \to \PrL$. 
Finally, one performs a Kan extension \cite[Proposition~A.5.12]{Mann:SixFunctors} to extend the functor from the non-full, wide subcategory $\Corr(\Sch^{\ft}_k,\textup{sep})$ to $\Corr(\Sch_k^\ft)$.

Recall from \thref{notation functors} that \eqref{equation:definition-SH}, in particular, encodes the existence of a closed symmetric monoidal structure on each category $\DA(X)$, adjunctions $(f^*,f_*)$ and $(f_!,f^!)$ for any map $f\colon Y\r X$ between finite type $k$-schemes satisfying proper base change and the projection formula.

As in \eqref{functor bx}, we have the exterior product functor
\eqn\label{functor bx SH}
-\bx-: \DA(X)\t_{\DA(k)}\DA(Y)\to \DA(X\x_kY),
\xeqn
and likewise for the ind-dualizable subcategories $\DA(-)^\indbl$ from \thref{definition indbl}.
Here and below we often drop the choice of coefficients $\Lambda$ from the notation.

\assu
\thlabel{assumption Lambda k}
Let $\Lambda\in \CAlg(\Spt)$ be connective. 
The exponential characteristic of $k$ is invertible in $\pi_0 \Lambda$ and $k$ has finite Galois cohomological dimension on discrete $\pi_0\Lambda$-modules.
\xassu

\exam\thlabel{example:cohomological-dimension}
The assumption is satisfied if $k$ is separably closed, or finite, or an unordered number field, or if $\pi_0 \Lambda$ is a $\Q$-algebra \cite[Remark~2.4.9(2)]{AyoubGallauerVezzani:RigidAnalyticMotives}.
\xexam

\prop\thlabel{DA compactly generated}
Under \thref{assumption Lambda k}, the functor $$\DA\colon \Corr(\Sch^\ft_k)\r \PrL$$ takes values in the \ii-category  $\PrLomega$ of compactly generated categories with functors preserving compact objects. Moreover, $\DA(k)\in \CAlg(\PrLomega)$ is locally rigid \cite[Definition~D.7.4.1]{Lurie:SAG}.
In particular, both $\DA$ and $\DA(-)^\indbl$ satisfy \thref{useful assumption}.
\xprop
\pf
Under \thref{assumption Lambda k} each category $\DA(X)$ is compactly generated \cite[Proposition~2.4.22(2), Remark~2.4.23]{AyoubGallauerVezzani:RigidAnalyticMotives}, the functors $f^*, f_!$ preserve compact objects and so does the exterior product \eqref{functor bx SH} by \cite[Theorem~2.4.9]{BondarkoDeglise:Homotopy-t-structure}.
In addition, $\DA(k)$ is locally rigid by \textit{loc.~cit.}, so we have $\DA(X)^{\textup{at}/\DA(k)}=\DA(X)^\omega$ (\thref{atomic vs compact}). 
This confirms \thref{useful assumption} for $\DA$, and also for $\DA(-)^\indbl$ by \thref{enriched maps ind dualizable}.
\xpf

The proof of \thref{DA compactly generated} shows that the atomic objects are precisely the compact ones.
These objects admit the following more geometric characterization (cf.~\cite[Theorem~6.3.26]{CisinskiDeglise:Etale} and \cite[Theorem~4.1]{RuimyTubach:EtaleMotives}):

\lemm\thlabel{characterization compact objects}
Let $X\in \Sch_k^\ft$. 
Under \thref{assumption Lambda k}, the following subcategories of $\DA(X)$ agree:
\begin{enumerate}
	\item the atomic objects $\DA(X)^{\textup{at}/\DA(k)}$,
	\item the compact objects $\DA(X)^\omega$,
	\item the constructible objects $\DA(X)^\cons$, i.e., those objects in $\DA(X)$ that on a locally closed decomposition of $X$ become dualizable, 
	\item the geometric objects in $\DA(X)^\gm$, i.e., the stable idempotent complete subcategory of $\DA(X)$ generated by the objects $f_\sharp 1_Y(n)$, where $f\colon Y \r X$ is smooth of finite type, $n \in \Z$.
\end{enumerate}
\xlemm
\pf
Using \thref{assumption Lambda k} we have $\DA(X)^{\textup{at}/\DA(k)}=\DA(X)^\omega=\DA(X)^\gm$: the first equality holds by the proof of \thref{DA compactly generated}; the second one by \cite[Proposition~2.4.22(2), Remark~2.4.23]{AyoubGallauerVezzani:RigidAnalyticMotives}.
So, it suffices to show the following inclusions:

$\DA(X)^\cons \subset \DA(X)^\omega$:
Let $M\in \DA(X)$ be constructible. 
Passing to a locally closed decomposition of $X$ witnessing the constructibility of $M$ and using the preservation of compact objects under $!$-push forward (\thref{DA compactly generated}) we may assume that $M$ is dualizable.
In this case its compactness follows from $\Hom_X(M,-)=\Hom_X(1_X,-\otimes M^\vee)$ and $1_X\in \DA(X)^\gm=\DA(X)^\omega$.

$\DA(X)^\omega \subset \DA(X)^\cons$: 
Let $M\in \DA(X)$ be compact.
Any residue field of $X$ is of finite transcendence degree over $k$, thus satisfies itself \thref{assumption Lambda k}.
We apply this to the generic point $\eta\r X$ of some irreducible component to conclude that $\DA(\eta)$ is locally rigid (\thref{DA compactly generated}).
In particular, the compact object $\eta^*M\in \DA(\eta)^\omega$ is dualizable.
Using continuity \cite[Theorem~2.5.1(2)]{AyoubGallauerVezzani:RigidAnalyticMotives} we find some open subset $U\subset X$ containing $\eta$ such that $M|_U$ is dualizable. 
So, $M$ is constructible by Noetherian induction.
\xpf

\subsection{A conjectural K\"unneth formula}
The exterior product functor is fully faithful on étale motives by the following result.

\prop\thlabel{proposition:SH-fully-faithful}
Under \thref{assumption Lambda k}, the exterior product \eqref{functor bx SH} for both $\DA$ and $\DA(-)^\indbl$ is fully faithful.
\xprop
\pf
\thref{DA compactly generated} confirms the categorical assumptions of \thref{proposition:abstract-fully-faithful} and \thref{proposition:abstract-fully-faithful indbl}; the key assumption on the isomorphism \eqref{proposition:fully-faithful:eq1} there is the content of \cite[Equations~(2.4.6.1) and (2.4.6.4)]{JinYang:Kuenneth}. 
\xpf

\rema\thlabel{fully faithful SH}
The consideration of the étale (as opposed to the Nisnevich) topology is not essential for \thref{proposition:SH-fully-faithful}. 
Indeed, for $\SH$ instead of $\DA$, the same arguments go through; in this case one may drop the condition that $k$ has finite cohomological dimension, see \cite[Théorème~4.5.67]{Ayoub:Six2} for the compact generation of $\SH(X)$. 
The citations \cite{BondarkoDeglise:Homotopy-t-structure,JinYang:Kuenneth} above only use the Nisnevich topology.
\xrema

As in \thref{definition indbl} we consider the subcategories $\DA(-)^{\indbl}$ of ind-dualizable objects.
Since $\DA(k)^\indbl=\DA(k)$ under \thref{assumption Lambda k} (\thref{DA compactly generated}) we can form the exterior product functor
\eqn\label{exterior product inddualizable}
\DA(X)^\indbl\t_{\DA(k)}\DA(Y)^\indbl\r \DA(X\x_kY)^\indbl,
\xeqn
which is fully faithful by \thref{proposition:SH-fully-faithful}.

\conj\thlabel{conjecture kuenneth}
Let $X, Y \in \Sch_k^\ft$ with $Y$ smooth over $k$.
Under \thref{assumption Lambda k}, the exterior product functor \eqref{exterior product inddualizable} is an equivalence in each of the following cases:
\begin{enumerate}
	\item \label{conjecture kuenneth 1} $\Char k=0$, or
	\item \label{conjecture kuenneth 2} $Y$ is proper over $k$.
\end{enumerate}
\xconj

Conditions \eqref{conjecture kuenneth 1} and \eqref{conjecture kuenneth 2} seem necessary in view of \thref{A2}.
In the following sections, we provide evidence for the essential surjectivity of \eqref{exterior product inddualizable}: in a nutshell, we prove in \thref{Künneth DMf indlis} that \thref{conjecture kuenneth} holds modulo the kernel of the adic realization (cf.~also \thref{Künneth Dad indlis} for the case of adic sheaves).
For finite fields $k$, the ind-dualizability can conjecturally be dropped by passing to Weil objects (\S\ref{sect--Künneth Weil motives}).

\exam
\thlabel{A2}
Let $k = \Fpq$, $X=Y=\A^1_k$ and consider the $\Z/p$-Galois cover 
\eqn
T:=\{t^p-t=xy\}\overset{f}{\to} X\times_{k} Y=\A^2_{k}, (x,y,t)\mapsto (x,y)
\xeqn
where $x,y$ are the coordinates on $X,Y$ respectively and $T\subset \A^3_{k}$ is the closed subscheme defined by displayed equation. 
Then, $f_\#\unit_T$ is dualizable, but does not lie in the essential image of $\DA(\A^1_k) \t_{\DA(k)} \DA(\A^1_k)$:
its mod-$n$ reduction for any $n\in \Z_{\geq 1}$ prime to $p$ does not lie in the essential image by \cite[Example~1.4]{HemoRicharzScholbach:Kuenneth} using rigidity (\thref{example:coefficients}\refit{Z/ell}).
\xexam

\section{The K\"unneth formula for adic sheaves}\label{section:adic-sheaves}
Throughout this section, let $k$ be a field of finite Galois cohomological dimension, so \thref{assumption Lambda k} holds for $\Lambda=\Z$. 
We consider the $6$-functor formalism for the category of étale motives
\begin{equation}\label{DA six functor}
\DA:=\DA_\et(-,\Z)\colon \Corr(\Sch_k^\ft)\r \PrLomegaZ,
\end{equation}
taking values in $\Z$-linear compactly generated presentable categories $\PrLomegaZ$, i.e., in modules over the presentably symmetric monoidal category $\Mod_\Z$ in compactly generated categories $\PrLomega$ (\thref{DA compactly generated}). 
Its exterior product is fully faithful (\thref{proposition:SH-fully-faithful}).
Note that $\DA(-)=\DA_\et(-,\Z[p^{-1}])$ if $\Char k= p>0$ by rigidity (\thref{example:coefficients}\eqref{Z/ell}). 

\subsection{The adic realization}\label{section:cohomological-motives-definition}

For each $n\in \Z_{\geq 1}$, the mod-$n$ functor $\PrLomegaZ\r \PrLomegaZ, \D\mapsto \D/n := \D \t_{\Mod_\Z} \Mod_{\Z/n}$ 
is a lax symmetric monoidal endofunctor.
Passing to the limit for varying $n$ we obtain a lax symmetric monoidal endofunctor
\begin{equation}\label{profinite completion}
\widehat{(-)}\colon \PrLomegaZ\r \PrLomegaZ, \;\; \D\mapsto \Dad
\end{equation}
where the limit in $\PrLomegaZ$ is computed as $\Dad=\Ind(\lim_{n, \Cat_\infty}\D^\omega/n)$.

\defi\thlabel{adic realization definition}
The \textit{formalism of adic sheaves} is the lax symmetric monoidal functor 
$$
\Dad:=\widehat\DA : \Corr(\Sch_k^\ft)\r \PrLomegaZ
$$
given by composing the formalism of étale motives \eqref{DA six functor} with the profinite completion \eqref{profinite completion}. 
The \textit{adic realization} is the induced transformation 
$$\DA\overset{\rho}{\r} \Dad$$ 
of lax symmetric monoidal functors $\Corr(\Sch_k^\ft)\r \PrLomegaZ$.
\xdefi

\rema
\thlabel{remark rho}
For $X\in \Sch_k^\ft$ the functor $\rho_X\colon \DA(X)\r \Dad(X)$ is characterized as the colimit preserving functor given by $\rho_X(M)= \lim_{n\geq 1} M/n$ on compact objects $M\in \DA(X)^\omega$.
\xrema

The adic realization $\rho$ is compatible with the $6$-functor formalism in the following sense.

\lemm\thlabel{remark:six-functors-realization}
The adic realization restricts to a natural transformation between lax symmetric monoidal functors
$$
\DA \overset{\rho}{\r} \Dad : \Corr(\Sch_k^\ft) \r \PrLomegaZ
$$
compatibly with the right adjoints $f_*$, $f^!$ for any map $f\colon Y\to X$ of finite type $k$-schemes and also $\IHom(-,-)$ when restricted to compact objects in the first variable.
In particular, one has natural equivalences 
\eqn\label{remark:six-functors-realization:eq1}
\begin{aligned}
& f^*\circ\rho_X=\rho_Y\circ f^*,\;\;\;  f^!\circ\rho_X=\rho_Y\circ f^!\\
\text{and}\;\;\; & f_*\circ\rho_Y=\rho_X\circ f_*,\;\;\; f_!\circ\rho_Y=\rho_X\circ f_!\\
\text{and}\;\;\; & \rho_X\circ \IHom_{\DA(X)}(-,-)=\IHom_{\Dad(X)}(-,-)\circ (\rho_X^\opp, \rho_X),
\end{aligned}
\xeqn
where one restricts to compact objects in the first variable of the inner homs.
\xlemm

\pf
This is a formal consequence of the definition of $\Dad(X)$ as the ind-completion of the profinite completion of $\DA(X)^\omega$ and the fact \cite[Théorèmes~8.10, 8.12]{Ayoub:Realisation} that the functors $f_*$ and $f^!$ and $\IHom$ (out of compact objects) preserve compact objects under our assumption on $k$.
\xpf

\rema\thlabel{formal adic realization}
More formally, \thref{remark:six-functors-realization} shows that the adic realization $\DA \overset{\rho}{\r} \Dad$ is an element in $\Fun'(\Corr(\Sch^\ft_k), \Fun(\Delta^1,\PrSt))$, using the notation around \eqref{Fun'}.
\xrema

Let us make the categories of adic sheaves more explicit. 
For each $X\in \Sch_k^\ft$ and each $n\in \Z$, rigidity (\thref{example:coefficients}\eqref{Z/ell}) implies that $\DA(X)/n=\HypShv(X_\et,\Z/n)$ if $\Char k\nmid n$ and $\DA(X)/n=0$ if $n$ is a $p$-power.
In particular, we see
\begin{equation}\label{DA rigidity}
\lim_{n} \DA(X)^\omega/n=\lim_{p\nmid n}\HypShv(X_\et,\Z/n)^\omega
\end{equation}
where the limits are formed in $\Cat_\infty$. 

Denote by $\Zhatp=\lim_{p\nmid n}\Z/n$ the prime-to-$p$ profinite completion of the integers viewed as a condensed ring.  
We consider the category $\HypShv(X_\proet,\Zhatp)$ of hypersheaves of $\Zhatp$-modules on the proétale site $X_\proet$ and its full subcategory $\Dcons(X,\Zhatp)$ of objects that are dualizable along a constructible stratification \cite{HemoRicharzScholbach:Constructible}. 
This defines a lax symmetric monoidal functor
\begin{equation}\label{ind constructible}
\Ind\big(\Dcons(-,\Zhatp)\big) \colon  \Corr(\Sch_k^\ft)\r \PrLomegaZ,
\end{equation}
which is constructed analogously to \eqref{equation:definition-SH} using \cite{BhattScholze:ProEtale}.

\lemm\thlabel{Dad comparison}
One has a natural equivalence 
\begin{equation}\label{indization constructible}
\Dad \cong \Ind\big(\Dcons(-,\Zhatp)\big) : \Corr(\Sch_k^\ft)\r \PrLomegaZ
\end{equation}
of lax symmetric monoidal functors.
In particular, for each $X\in \Sch_k^\ft$, the category $\Dad(X)$ is the full subcategory of $\HypShv(X_\proet,\Zhatp)$ generated under filtered colimits by constructible objects, and the natural $t$-structure on $\HypShv(X_\proet,\Zhatp)$ restricts to $\Dad(X)$.
\xlemm
\pf
One has natural lax symmetric monoidal transformations 
\[
\Dad^\omega \to \lim_{p\nmid n}\HypShv((-)_\et,\Z/n)^\omega \r \lim_{p\nmid n}\Dcons(-,\Z/n) \gets \Dcons(-,\Zhatp),
\]
where the first is an equivalence by rigidity \eqref{DA rigidity}, the second by \cite[Proposition~7.1]{HemoRicharzScholbach:Constructible} and \cite[Proposition~6.4.8]{BhattScholze:ProEtale} (using that $k$ has finite Galois cohomological dimension), and the third by \cite[Proposition~5.1]{HemoRicharzScholbach:Constructible}.
Passing to ind-completions yields the equivalence \eqref{indization constructible}.
Thus, $\Dad(X)$ defines a full subcategory of $\HypShv(X_\proet,\Zhatp)$ by \cite[Corollary~8.3]{HemoRicharzScholbach:Constructible}, again using that $k$ has finite Galois cohomological dimension. 
Since $\Zhatp$ is t-admissible \cite[Corollary~6.11(2)]{HemoRicharzScholbach:Constructible}, the t-structure on $\HypShv(X_\proet,\Zhatp)$ restricts to $\Dad(X)^\omega=\Dcons(X,\Zhatp)$ by \cite[Theorem~6.2(3)]{HemoRicharzScholbach:Constructible} and thus also to $\Dad(X)$ by t-exactness of filtered colimits.
The compatibility of the equivalence with $f^*$, $f_!$ and $\bx$ holds by \cite[Théorème~6.6]{Ayoub:Realisation}.
\xpf

\rema
\thref{Dad comparison} asserts that the functors $f^*$, $f_!$ and $\bx$ on $\Dad$ induced by the ones on $\DA$ are the same functors as the ones on $\Ind\big(\Dcons(-,\Zhatp)\big)$ induced from the sheaf-theoretic functors for étale torsion sheaves.
Since it is an equivalence, this is also true for their right adjoints.
\xrema

\subsection{Ind-dualizable adic sheaves}
As in \thref{definition indbl}, we consider the subcategories $\Dad(-)^{\indbl}$ of ind-dualizable objects.
Under the equivalence in \thref{Dad comparison}, this is the category of ind-lisse sheaves $\D_{\textrm{indlis}}(-,\Zhatp)$ introduced in \cite{HemoRicharzScholbach:Constructible}. 

\prop
\thlabel{Dad pullback conservative}
For $X\in \Sch_k^\ft$ connected and a rational point $x: \Spec k \r X$, the functor
$$(x \x \id_Y)^*\colon \Dad(X \x_k Y)^\indbl \r \Dad(Y)$$
is conservative.
In particular, $(x \x \id_Y)^!$ is conservative if, additionally, $x\r X$ lies in the regular locus.
\xprop

\pf
The final statement follows from the first statement by continuity of $!$-pullback and purity. 
For the first statement, we use t-structures as follows.

The natural $t$-structure on hypersheaves restricts to the categories $\Dad(-)^\indbl$ because it does so on dualizable objects \cite[Theorems 6.2, 6.12]{HemoRicharzScholbach:Constructible} (by t-admissibility of $\Zhatp$ \cite[Lemmas~6.7(3), 6.10]{HemoRicharzScholbach:Constructible}) and because filtered colimits are t-exact.
(As in \thref{Dad comparison} the category $\Dad(X)^\indbl$ is compactly generated by the dualizable objects $\Dad(X)^\dbl$.)
Since the natural t-structure on hypersheaves is non-degenerate and $*$-pullback is t-exact, it is enough to show conservativity for the hearts $\Dad(-)^{\indbl,\heartsuit}=\Ind(\Dad(-)^{\indbl,\heartsuit})$.

So, let $M\in \Dad(X)^{\indbl,\heartsuit}$ such that $x^*M=0$.
If $M$ is dualizable, the local constancy of dualizable sheaves \cite[Theorem~4.13]{HemoRicharzScholbach:Constructible} implies $M=0$ using the connectedness of $X$.
The general ind-dualizable case now follows from \thref{lemma-ind-conservativity}.
\xpf

\lemm\thlabel{lemma-ind-conservativity}
Let $F\colon A\r B$ be a conservative, exact functor between abelian categories. 
Then, its indization $\Ind F\colon \Ind A\r \Ind B$ is conservative. 
\xlemm
\pf
An exact functor between abelian categories is conservative if and only if it is faithful. 
So, it suffices to show that $\tilde F:=\Ind F$ is faithful if $F$ is so. 
Given a filtered colimit $M=\colim M_i$ with $M_i\in A$ and $N\in A$, we get the injection
\[
\colim \Hom_A(N,M_i)\r \colim \Hom_{\Ind B}(F(N), F(M_i)),
\]
by exactness of filtered colimits in abelian groups. 
The claim for general $N\in \Ind A$ follows formally because limits preserve injections of abelian groups. 
\xpf

\prop
\thlabel{pushforward char 0}
Suppose $\Char k=0$.
Let $X, Y \in \Sch_k^\ft$ with $Y$ smooth over $k$.
Then, the functor
$$\pr_{X,*}\colon \Dad(X \x Y) \r \Dad(X)$$
preserves dualizable objects.
\xprop

\pf
Let $M\in \Dad(X\x Y)$ be a dualizable object. 
Then, $M$ is constructible and so is $\pr_{X,*}M$.
In order to show that $\pr_{X,*}M$ is dualizable we reduce to the case $M=\Z/\ell$ for some prime $\ell\neq p$ as follows. 

The functor
\[
\Dad_\cons(X) \to \Dad_\cons(X)/\ell=\D_\cons(X,\Z/\ell),\pr_{X,*}M\mapsto \pr_{X,*}M/\ell
\]
is symmetric monoidal, commutes with inner homs (because one has $\IHom(-,-)/\ell=\IHom(-,-/\ell)$ as holds for any exact functor; then use standard $\otimes$-$\Hom$-adjunctions) and the family for all primes $\ell\neq p$ is jointly conservative (because constructible objects are profinitely derived complete).
Thus, we may assume $M=M/\ell$ for some $\ell\neq p$.
Since dualizability is proétale local \cite[Lemma~4.5]{HemoRicharzScholbach:Constructible} we may assume that $k$ is algebraically closed and $X, Y$ are connected.
Using the t-structure on dualizable $\bbZ/\ell$-sheaves \cite[Theorem~6.2(2)]{HemoRicharzScholbach:Constructible} and an induction on the amplitude we reduce to the case that $M$ is concentrated in a single cohomological degree and after shifting in degree $0$, i.e., $M$ is an étale locally constant sheaf of finite dimensional $\bbZ/\ell$-vector spaces \cite[Proposition~7.1]{HemoRicharzScholbach:Constructible}.  
Fixing some geometric point in $X\x Y$, the sheaf $M$ corresponds under the monodromy equivalence to a continuous representation of the étale fundamental group
\[
\pi_1^\et(X\x Y)\overset{\cong}{\lr} \pi_1^\et(X)\x \pi_1^\et(Y),
\]
where the K\"unneth formula holds because $\Char k=0$ (via comparison with the topological fundamental group \cite[Exp.~XII, Corollaire~5.2]{SGA1}; alternatively see \cite[Corollary~2.12]{Liu:Kuenneth} for an algebraic argument).
Since $\bbZ/\ell$ is finite there are finite étale morphisms $\tilde X\r X$, $\tilde Y\r Y$ with $\tilde X$, $\tilde Y$ connected such that $M|_{\tilde X\x \tilde Y}$ is constant.
Then, we reduce to the case $\tilde Y=Y$ by étale descent and to the case $\tilde X=X$ by base change.
That is, we may assume that $M$ is a constant $\Z/\ell$-vector space and further that $M$ is the constant sheaf $\bbZ/\ell$ on $X\x Y$ sitting in degree $0$.
This finishes the reduction step.

Further, the formation of $\pr_{X,*}\Z/\ell$ commutes with base change (i.e., for $g \colon X' \r X$, $\pr_{X'} \colon Y \x_k X' \r X'$, we have $g^* \pr_{X,*} \Z/\ell = \pr_{X',*} \Z/\ell$; in other words, the pair $(\pr_X, \Z/\ell)$ is cohomologically proper in the sense of \cite[Th.~finitude, Appendix]{SGA412}).
Indeed, this follows, for example, from factoring the smooth map $Y\to \Spec k$ as a chain of elementary fibrations \cite[Exp.~XI, 3.3]{SGA4:3} 
and using the hypothesis $\Char k = 0$ to see that there is only tame ramification in order to apply \cite[Th.~finitude, Appendix, 5.1.3]{SGA412}.

Finally, for each étale specialization $\bar t\leadsto \bar s$ of geometric points in $X$ there is the commutative diagram
$$\xymatrix{
(\pr_{X,*}\bbZ/\ell)_{\bar s}\ar[d]^{\cong} \ar[r]^{\textrm{specialization}} & (\pr_{X,*}\bbZ/\ell)_{\bar t} \ar[d]^{\cong} \\
\RG(\bar s\x Y,\bbZ/\ell) \ar[r]^{\cong} & \RG(\bar t\x Y,\bbZ/\ell),
}$$
where the vertical maps are isomorphisms by base change and the bottom map by independence of mod-$\ell$ étale cohomology on the choice of an algebraically closed overfield of $k$.
Hence, the constructible sheaf $\pr_{X,*}\bbZ/\ell$ is étale locally constant (hence, dualizable) by \cite[\href{https://stacks.math.columbia.edu/tag/0GJ7}{0GJ7}]{StacksProject} applied to the cohomology sheaves.  
\xpf

\subsection{The Künneth formula for ind-lisse sheaves}
The following theorem confirms the analogue of \thref{conjecture kuenneth} for adic sheaves.
A similar result also holds for all primes $\ell \ne \text{char}(k)$ in the $\ell$-adic formalism $\Ind\Dcons(-,\Zl)$ obtained by $\ell$-completing $\Dad$ in $\PrLomegaZ$.

\theo
\thlabel{Künneth Dad indlis}
Let $k$ be a field of finite Galois cohomological dimension and $X, Y \in \Sch^\ft_k$ with $Y$ smooth over $k$.
Assume one of the following:
\begin{enumerate}
    \item \label{char zero} $\Char k = 0$, or
    \item \label{proper} $Y$ is proper over $k$.
\end{enumerate}
Then, the exterior product functor
\begin{equation}\label{Künneth Dad indlis eq}
\Dad(X)^\indbl \t_{\Dad(k)} \Dad(Y)^\indbl  \r \Dad(X \x_k Y)^\indbl
\end{equation}
is an equivalence.
\xtheo
\pf
The functor \eqref{Künneth Dad indlis eq} is fully faithful (also on all adic sheaves $\Dad$) by the proof of \thref{proposition:SH-fully-faithful}.
The categorical assumptions hold by construction \eqref{indization constructible}.

For essential surjectivity of \eqref{Künneth Dad indlis eq} we may assume that $X$ is connected. 
We first assume that there is a rational point $x \colon \Spec k \r X$ supported in the regular locus.
In this case we apply Propositions \ref{proposition:essentially-surjective} and \ref{Künneth inddualizable} to $\Dad\colon \Corr(\Sch^\ft_k) \r \PrgmL$ and check the conditions in these propositions:
\begin{itemize}
\item\textit{\thref{proposition:essentially-surjective}\eqref{pullback conservative}: $(x \x \id_Y)^* : \Dad(X \x_k Y)^\indbl \r \Dad(Y)$ is conservative.}
This is \thref{Dad pullback conservative}.

\item\textit{\thref{Künneth inddualizable}\eqref{pr* dualizable}: $\pr_{X,*}\colon \Dad(X \x_k Y) \r \Dad(X)$ preserves dualizable objects.}
If $\Char k=0$, this holds by \thref{pushforward char 0}. 
If $Y$ is proper over $k$, then $\pr_{X}$ is smooth and proper, say of relative dimension $d$, this follows from relative purity: for dualizable $M$, $\pr_{X,*} M$ is dualizable with dual $$\pr_{X,*} (M^\dual)(d)[2d] = \pr_{X,*} \IHom_{X\x Y}(M, \pr_X^! 1_X) = \IHom_X(\pr_{X,!} M, 1_X).$$

\item\textit{\thref{Künneth inddualizable}, Items \eqref{x pure} and \eqref{*! exchange}: $x^!1_X$ is $\t$-invertible and $ p_Y^* x^! 1_X \cong (x \x \id_Y)^! \pr_X^* 1_X$.} 
This holds by absolute purity $x^! 1_X \cong 1_x(-d_x)[-2 d_x]$ (with $d_x=\dim_{\textup{Krull}}\calO_{X,x}$) and its compatibility under transversal pullback (\thref{purity remark}).
\end{itemize}

We now handle the case where the regular locus of $X$ does not have a rational point. 
Then, there is a finite separable field extension $K\supset k$ such that the regular locus of $(X \x_k K)_\red$ admits a $K$-rational point \StP{056U}.
Note that the Künneth formula 
\eqn\label{kuenneth for field extensions}
\Dad(Z)^{(\indbl)} \t_{\Dad(k)} \Dad(K) = \Dad(Z \x_k K)^{(\indbl)} 
\xeqn
holds for all $Z\in \Sch_k^\ft$ by \thref{pushforward conservative} because the pushforward $(\id_Z \x (\Spec K \r \Spec k))_*$ preserves dualizable objects and is conservative.
(Conservativity is a general property of a $6$-functor formalism for which pullback along étale surjective maps is conservative, cf.~the proof of \cite[Proposition~2.4.6]{Cisinski:Cohomological}.)
We combine \eqref{kuenneth for field extensions} with the previous step (all undecorated tensor products are over $\Dad(k)$):
$$\xymatrix{\Dad(K) \t \Dad(X)^{\indbl} \t \Dad(Y)^{\indbl} \ar[d]^\cong \\
\Dad(X_K)^\indbl \t_{\Dad(K)} \Dad(Y_K)^\indbl \ar[r]^{\cong\phantom{hhhhh}} & \Dad(X_K \x_K Y_K)^\indbl = \Dad(X \x_k Y \x_k K)}$$
The Künneth formula for $X, Y$ follows from the one for the triplet $X, Y, \Spec K$ by the argument of \thref{Kuenneth constructible}, given that the pullback functor $\Dad(X \x_k Y) \r \Dad(X \x_k Y \x_k K)$ is conservative.
\xpf

\section{The K\"unneth formula for cohomological motives}\label{section:cohomological-motives}
As in \S\ref{section:adic-sheaves} we continue to fix a field $k$ of finite Galois cohomological dimension and consider the $6$-functor formalism for the category of étale motives with $\Z$-coefficients \eqref{DA six functor}:
\begin{equation}\label{DA six functors 2}
\DA:=\DA_\et(-,\Z)\colon \Corr(\Sch_k^\ft)\r \PrLomegaZ
\end{equation}
The following definition is analogous to \cite[Definition~2.6]{Ayoub:Motives} in the Betti setting.

\defi\thlabel{definition:cohomological-motives}
For $X\in \Sch_k^\ft$, the \ii-category of \textit{cohomological motives} is the Verdier quotient (cf.~\S\ref{sect--enforcing conservativity}) in $\PrSt$,
\eqn
\DAf(X):=\DA(X)/\ker \rho_X,
\xeqn
where $\rho_X\colon \DA(X) \r \Dad(X)$ is the adic realization (\thref{adic realization definition}) and $\ker \rho_X$ is the full subcategory of objects $M\in\DA(X)$ with $\rho_X(M)= 0$.
\xdefi 

\rema
The homotopy category $\Ho\DAf(X)$ is triangulated and identifies with the Verdier quotient $\Ho\DA(X) / \Ho\ker(\rho_X)$ of triangulated categories \cite[Proposition~5.9]{BlumbergGepnerTabuada:UniversalKTheory}.
\xrema

By design (cf.~\S\ref{sect--enforcing conservativity}), the category of cohomological motives fits into a triangle
\begin{equation}
\xymatrix{
\DA(X) \ar[dr]_\flat \ar[rr]^{\rho_X} & & \Dad(X) \\
& \DAf(X) \ar[ur]_{\rho_X^\flat}
}
\label{diagram DAX}
\end{equation}
where the functor $\flat$ is a localization and $\rho_X^\flat$ is conservative.
In \cite[Conjecture~2.12]{Ayoub:Motives}, Ayoub raises the question whether the Betti realization (for $k$ a field of characteristic 0 with finite transcendence degree over $\mathbb Q$) is conservative on the entire category $\DA(X)$, which is even stronger than asking it to be conservative on the subcategory $\DA(X)^\comp$ of compact objects.
The functor $\flat$ is an equivalence precisely if the analogous statement holds for the adic realization.
Without addressing these questions, let us point out that the formation of $\DAf(X)$ only concerns rational (and not modular) coefficients.

\lemm\thlabel{lemma:detect-equivalences}
Let $X\in \Sch_k^\ft$ and $M\in \DA(X)$.
Then, the following are equivalent:
\begin{enumerate}
	\item \label{lemma:detect-equivalences:it1}
			$M^\flat=0$;
	\item \label{lemma:detect-equivalences:it2}
			$\rho_X(M)=0$;
	\item \label{lemma:detect-equivalences:it3}
			$\rho_X(M_\Q)=0$ and $M/\ell=0$ for all primes $\ell\neq \Char k$.
\end{enumerate}
In addition, if $M\in \DA(X)^\omega$, then \eqref{lemma:detect-equivalences:it1}--\eqref{lemma:detect-equivalences:it3} are equivalent to $M/\ell=0$ for all primes $\ell\neq p$.
\xlemm
\pf
\eqref{lemma:detect-equivalences:it1}$\iff$\eqref{lemma:detect-equivalences:it2}: This holds by the universal property of Verdier quotients.

\eqref{lemma:detect-equivalences:it2}$\iff$\eqref{lemma:detect-equivalences:it3}:
Note that $M/\ell=0$ for all primes $\ell\neq p$ is equivalent to $M/n=0$ for all $0\neq n\in \Z$ using rigidity (\thref{example:coefficients}\eqref{Z/ell}) if $\Char k=p>0$ and $n$ is a $p$-power.
Further, one has $M/n=0$ if and only if $\rho_X(M/n)=0$ because $\rho_X/n$ is an equivalence (\thref{remark rho}). 
In particular, $\ker\rho_X \subset \D(X)_\Q$, which shows the equivalence of Items \eqref{lemma:detect-equivalences:it2} and \eqref{lemma:detect-equivalences:it3}. 

Finally, if $M\in \DA(X)^\omega$ and $M/\ell=0$ for all $\ell\neq \Char k$, then $M/n=0$ for all $0\neq n\in \Z$ and so $\rho_X(M)=\lim_{n\in \Z}M/n=0$ (\thref{remark rho}).

\xpf

\rema
\thlabel{omega_1 compact}
The adic realization $\rho_X$ is a functor in $\PrStOm$ (\thref{adic realization definition}).
Hence, $\ker\rho_X$ is $\omega_1$-compactly generated and objects $M\in (\ker \rho_X)^{\omega_1}$ can be written as sequential colimits
\eqn
M=\colim (M_1\overset{f_1}{\lr}M_2\overset{f_2}{\lr}\cdots)
\xeqn
with $M_i\in \DA(X)^\omega$ and $\rho(f_i)=0$ for all $i\geq 1$  \cite[Theorem~7.4.1]{Krause:Localization}.
In particular, $\DAf(X)$ is $\omega_1$-compactly generated but we do not know whether it is compactly generated. 
On the bright side, we will see in \thref{DM gm flat} below that $\DAf(X)$ is $\DAf(k)$-atomically generated, which suffices to prove the categorical K\"unneth formula.
\xrema

\subsection{The six functor formalism}\label{section:six-functors}
Recall the characterization of geometric objects $\DA(X)^\gm$ from \thref{characterization compact objects}.
Let us denote by $\DAf(X)^\gm$ the essential image of the composition $\DA(X)^\gm\subset \DA(X)\to \DAf(X)$. 

\prop\thlabel{six-functors-SH-flat}
The assignment $X \mapsto \DAf(X)$ is part of a lax symmetric monoidal functor
$$
\DAf \colon \Corr(\Sch^\ft_k) \r \PrL.
$$
There is a natural transformations of lax symmetric monoidal functors
$$\DA \stackrel \flat \r \DAf \stackrel{\rho^\flat} \r \Dad$$
factorizing the adic realization $\rho$ (\thref{adic realization definition}).
The natural transformations $\flat,\rho^\flat$ are compatible with the right adjoints $f_*$, $f^!$ and also $\IHom$ when restricted in the first variable to objects in $\DA^\gm$ and $(\DAf)^\gm$ respectively. 
\xprop

\pf
Consider $F := \Fun(\Corr(\Sch_k^\ft), \PrLSt)$, equipped with its Day convolution product. Commutative algebra objects in $F$ are precisely lax symmetric monoidal functors such as $\DA$ or $\Dad$, and $\rho$ is an object in $\Fun(\Delta^1, F)$.
We apply the method of \thref{factorization}  (i.e., replace $\PrLSt$ there by $F$). Using that kernels and pushouts in $F$ are computed pointwise, we obtain a map $\Fun(\Delta^1, F) \r \Fun(\Delta^2, F)$
that sends $\rho$ to a functor whose evaluation at some object $X \in \Sch_k^\ft$ is the diagram \eqref{diagram DAX}.

The compatibility with $f_*, f^!$ follows from \thref{factorization and right adjoints} and \thref{formal adic realization}, i.e., the adic realization functor $\rho\colon \DA \r \Dad$ is an element in $\Fun'(\Corr(\Sch^\ft_k), \Fun(\Delta^1,\PrSt))$, using the notation around \eqref{Fun'}.
The compatibility of $\flat$ with $\IHom$ as stated follows from the same considerations: in \eqref{vertically right adjointable} consider the situation where $\sigma_0$ is given by $-\t_X M\colon \DA(X)\r\DA(X)$ for some $M\in \DA(X)^\gm$.
The functor preserves geometric (=compact by \thref{characterization compact objects}) objects \cite[Théorème~8.12]{Ayoub:Realisation}, so its right adjoint $\IHom_{\DA(X)}(M,-)$ is continuous and commutes with the adic realization (\thref{remark:six-functors-realization}).
This formally implies $$\IHom_X(M,-)^\flat=\IHom_X(M^\flat,(-)^\flat),\;\; \rho_X^\flat\circ \IHom_X(M^\flat,-)=\IHom_X(\rho_X(M), \rho_X^\flat(-))$$
as functors $\DA(X)\r \DAf(X)$ and $\DAf(X)\r\Dad(X)$ respectively.
\xpf

The proof, in particular, shows that the right adjoint functors $f_*, f^!, \IHom(M,-)$ (for $M$ geometric) when formed in the formalism $\DAf$ are continuous and induced from the universal property of the Verdier quotient by the kernel of the adic realization.    
These functors are even internally left adjoint by the following result.

\lemm
\thlabel{DM gm flat}
The formalisms $\DA^\flat$ and $\DA^\flat(-)^\indbl$ satisfy \thref{useful assumption}, i.e., take values in atomically generated $\DAf(k)$-modules.
\xlemm
\pf
The category $\DAf(X)$ is generated under colimits by $\DAf(X)^\gm$, as the same is true for $\DA$ in place of $\DAf$ and $(-)^\flat$ commutes with colimits.
Moreover, all objects in $\DA(k)^\gm=\DA(k)^\omega$ are dualizable (\thref{DA compactly generated}) and so are their images in $\DAf(k)$ by monoidality of $(-)^\flat$, i.e., $\DA^\flat(k)^\indbl=\DAf(k)$.
Applying \thref{atomic vs compact}\eqref{atomic vs compact 1} to $\DAf(k)\in \CAlg(\PrL)$ and observing that $p_{X,*}\IHom_X(M, -)$ is continuous ($M$ geometric) shows that all objects in $\DAf(X)^\gm$ are atomic.
Finally, since the functors $f^*, f_!$ and the exterior product on $\DA$ preserve geometric objects, the same is true for these functors on $\DAf$ by \thref{six-functors-SH-flat}, thus they preserve atomic objects (\thref{internal left adjoint lemma}).
The case of $\DAf(-)^\indbl$ follows from \thref{enriched maps ind dualizable}.
\xpf

\rema\thlabel{DAf dualizable}
In light of \thref{dualizability ramzi}, \thref{DM gm flat} shows that $\DAf(X)$ and $\DAf(X)^\indbl$ is dualizable in $\Mod_{\DAf(k)}(\PrL)$ but we do not know whether the same is true in $\PrSt$.
\xrema

We now list (only) those properties of $\DAf$ that we need in the sequel.
In fact, most of the properties for motives as surveyed, say, in \cite[Synopsis~2.1.1]{RicharzScholbach:Intersection} hold for $\DAf$ as well: properties i)--v) (six functors and invertibility of $\Z(1)$), vii)--xii) (projection formula, Verdier duality, homotopy invariance, purity), xiv) and xvi) (h-descent, realization functors) formulated there hold for $\DAf$ as well.
\thref{DM gm flat} replaces property vi) there.
We do not know how to compute, say, $\pi_0\End_{\DAf(k)}(1)$ (cf.~property xiii)), nor do we have a weight structure (cf.~property xv)).

\prop
\thlabel{DMf localization}
Cohomological motives satisfy the following properties:
\begin{enumerate}
    \item \label{DMf localization 1}
    (Localization) For a closed immersion $i\colon Z\r X$ in $\Sch_k^\ft$ with complementary open immersion $j\colon U\r X$, there are the following functorial fiber sequences in $\DAf(X)$:
$$j_! j^! \r \id \r i_* i^*,$$
$$i_! i^! \r \id \r j_* j^*.$$
\item\label{relative purity eq}
(Relative purity)
For a smooth map $f\colon Y\r X $ in $\Sch_k^\ft$ of relative dimension $d$, the natural transformation
\begin{equation}
f^*(-)\t f^!1_X\r f^!\colon \DAf(Y) \r \DAf(X). \label{relative purity}    
\end{equation}
is an isomorphism.
In addition, there is a canonical Thom isomorphism $f^!1_X\cong 1_Y[2d](d)$.
\item\label{absolute purity eq}
(Absolute purity)
For a regular immersion $i\colon Z\r X$ in $\Sch_k^\ft$ of codimension $d$, the family of motivic cohomology classes in \cite[Proposition~1.3.21]{BondarkoDeglise:Homotopy-t-structure} for $\DA$ induces an isomorphism
\eqn
1_Z[-2d](-d)\r i^!1_X \colon \DAf(X)\r \DAf(Z)
\xeqn
that is compatible with transversal pullback.  
\end{enumerate}
\xprop

\pf
This follows from the corresponding property of $\DA$ (cf., e.g., \cite[Theorem~2.4.50]{CisinskiDeglise:Triangulated} for localization and relative purity; \cite[Proposition~1.3.21]{BondarkoDeglise:Homotopy-t-structure} for absolute purity) and the fact that $\flat$ commutes with these functors.
\xpf

\coro\thlabel{atomic vs constructible}
For $X\in \Sch_k^\ft$, the category $\DAf(X)$ is $\DAf(k)$-atomically generated by its constructible objects $\DAf(X)^\cons$, i.e., those objects in $\DAf(X)$ that on a locally closed decomposition of $X$ become dualizable.
\xcoro
\pf
Any dualizable object is atomic (Lemmas~\ref{enriched maps ind dualizable} and \ref{DM gm flat})
and so is any constructible object by localization (\thref{DMf localization}\eqref{DMf localization 1}), given that the functor $\iota_!\iota^*\colon \DAf(X)\r\DAf(X)$ preserves atomic objects (\thref{DM gm flat}) for any locally closed immersion $\iota\colon Z\r X$.
This shows $\DAf(X)^\cons\subset \DAf(X)^{\textup{at}/\DAf(k)}$.
Since $\DA(X)^\cons$ generates $\DA(X)$ under colimits (\thref{characterization compact objects}) and $\flat$ preserves constructible objects, the same is true for $\DAf(X)^\cons\subset \DAf(X)$.
\xpf

\subsection{The Künneth formula for ind-dualizable cohomological motives}
We have now all tools for the proof of the K\"unneth formula.

\prop\thlabel{fully faithful DAf}
For $X, Y \in \Sch^\ft_k$, the exterior product functor
$$\DAf(X) \t_{\DAf(k)} \DAf(Y) \r \DAf(X \x_k Y)$$
is fully faithful, and likewise for ind-dualizable categories $\DAf(-)^\indbl$.
\xprop

\pf
This is analogous to the proof of \thref{proposition:SH-fully-faithful}:
\thref{DM gm flat} confirms the categorical assumptions of \thref{proposition:abstract-fully-faithful} and \thref{proposition:abstract-fully-faithful indbl}; the key assumption on the isomorphism \eqref{proposition:fully-faithful:eq1} there is the content of \cite[Equations~(2.4.6.1) and (2.4.6.4)]{JinYang:Kuenneth} for $\DA$ and carries over to $\DAf$ by \thref{six-functors-SH-flat}. 
\xpf

The next result confirms \thref{conjecture kuenneth} modulo the kernel of the adic realization:

\theo
\thlabel{Künneth DMf indlis}
Let $X, Y \in \Sch^\ft_k$ with $Y$ smooth over $k$.
Assume one of the following:
\begin{enumerate}
    \item \label{char zero} $\Char k = 0$, or
    \item \label{proper} $Y$ is proper over $k$.
\end{enumerate}
Then, the exterior product functor
\begin{equation}
\DAf(X)^\indbl \t_{\DAf(k)} \DAf(Y)^\indbl \r \DAf(X \x_k Y)^\indbl	
\label{Künneth DMf indlis eqn}
\end{equation}
is an equivalence.
\xtheo
\pf
The proof of \thref{Künneth Dad indlis} applies: 
First off, \eqref{Künneth DMf indlis eqn} is fully faithful by \thref{fully faithful DAf}.
Essential surjectivity in the case where the regular locus of $X$ admits a $k$-rational point follows from applying \thref{proposition:essentially-surjective} and \thref{Künneth inddualizable} to $\DAf$ and its ind-dualizable objects. 
For convenience let us swiftly check the conditions from these propositions:
\begin{itemize}
\item\textit{\thref{proposition:essentially-surjective}\eqref{pullback conservative}: $(x \x \id_Y)^* : \DAf(X \x_k Y)^\indbl \r \DAf(Y)$ is conservative ($X$ connected).}
This follows from the case of adic sheaves $\Dad$ (\thref{Dad pullback conservative}), the compatibility of pullbacks under the adic realization $\rho^\flat$ (\thref{six-functors-SH-flat}) and its conservativity, which holds by design of cohomological motives.

\item\textit{\thref{Künneth inddualizable}\eqref{pr* dualizable}: $\pr_{X,*}\colon \DAf(X \x_k Y) \r \DAf(X)$ preserves dualizable objects.}
If $\Char k=0$, this holds by \thref{pushforward char 0}, noting that $\rho^\flat$ commutes with $\pr_{X,*}$ and reflects whether an atomic objects is dualizable (because it is monoidal, preserves inner homs and is conservative). 
If $Y$ is proper over $k$, then $\pr_{X}$ is smooth and proper.
This case formally follows from relative purity as in \thref{DMf localization}\eqref{relative purity eq} (cf.~the proof of \thref{Künneth Dad indlis}).

\item \textit{\thref{Künneth inddualizable}\eqref{x pure} and \eqref{*! exchange}: $x^!1_X$ is $\t$-invertible and $ p_Y^* x^! 1_X \cong (x \x \id_Y)^! \pr_X^* 1_X$.} 
This holds by absolute purity and its compatibility under transversal pullback (\thref{DMf localization}\eqref{absolute purity eq}).
\end{itemize}

The general case follows as in the proof of \thref{Künneth Dad indlis} from a Galois descent argument.
\xpf 

\rema
\thlabel{Künneth formula prestacks}
\thref{DAf dualizable} implies that both sides in \eqref{Künneth DMf indlis eqn} turn colimits of prestacks (i.e., colimits in $\Fun((\Sch_k^\ft)^\opp, \Ani)$) into limits in $\Mod_{\DAf(k)}(\PrL)$. 
Therefore, the statement above extends immediately to arbitrary prestacks $X$ and prestacks $Y$ that admit a presentation of the form $Y = \colim Y_i$, where $Y_i \r \Spec k$ is smooth of finite type (resp.~smooth and proper if $\Char k > 0$).
\xrema

\section{The Künneth formula for Weil objects}
\label{sect--Künneth Weil motives}
Throughout this section, let $k=\Fq$ be a finite field of characteristic $p>0$ and cardinality $q$. We fix an algebraic closure $\bar k$.
In Section~\ref{Weil formalism} we review how to construct a formalism of Weil objects based on a $6$-functor formalism on $k$-schemes.
We prove the K\"unneth formula for adic Weil sheaves and for cohomological Weil motives (\thref{Künneth DMf Weil}), leading to \thref{conjecture kuenneth weil} for étale Weil motives.

\subsection{Weil formalism}\label{Weil formalism}
We consider the category $\PreStk_\bk:=\Fun((\Sch_\bk^\ft)^\opp, \Ani)$ of prestacks, which contains $\Sch_\bk^\ft$ as a full subcategory via the Yoneda embedding.
For $X\in \Sch_k^\ft$ the associated \textit{Weil prestack} is the $\bk$-prestack
$$X^W := \colim \left (\ol X \stackrel[\id]{\varphi_X} \rightrightarrows \ol X \right)\in \PreStk_{\bk},$$
where $\ol X:=X\x_k\bk$ and $\varphi_X:=\Frob_X \x_k \id_{\bk}$ with $\Frob_X$ the $q$-Frobenius on $X$.
The assignment $X\mapsto X^W$ induces a functor 
\eqn
(-)^W\colon \Sch_k^\ft\r \PreStk_\bk,
\xeqn
which is neither full, nor faithful, nor does it preserve products.
The functor formally extends to a symmetric monoidal faithful functor
\eqn\label{weil symmetric}
(-)^W\colon (\Sch_k^\ft)^{\sqcap}\r \PreStk_\bk, (X_i)_{i\in I}\mapsto \prod_{i\in I}X_i^W, 
\xeqn
where $(\Sch_k^\ft)^\sqcap$ is the free completion of $\Sch_k^\ft$ under finite products: its objects are finite families $(X_i)_{i\in I}$ and its morphisms $\Hom((X_i)_{i}, (Y_j)_{j}) = \prod_j \coprod_i \Hom(X_i, Y_j)$.
The category $(\Sch_k^\ft)^\sqcap$ has pullbacks and a terminal object (namely the empty family for the empty product is a singleton), thus finite limits.
The functor $\Sch_k^\ft\r (\Sch_k^\ft)^\sqcap, X\mapsto (X)$ preserves pullbacks, but \emph{not} finite products.

\lemm
Any lax symmetric monoidal functor $\D \colon (\Sch_\bk^\ft)^\opp \r \PrL$ induces a lax symmetric monoidal functor $\tD((-)^W)\colon (\Sch_k^\ft)^{\sqcap,\opp}\r \PrL$ satisfying 
\eqn\label{DXW}
\tD(X^W) = \lim \left(\D(\ol X) \stackrel[\id]{\varphi_X^*} \rightrightarrows \D(\ol X)\right)
\xeqn
for any $X\in \Sch_k^\ft$.
\xlemm

\pf
The right Kan extension $\tD\colon \PreStk_{\bk}^\opp\r \PrL$ of $\D$ to prestacks is lax symmetric monoidal by \cite[Theorem~8.5.3]{Hinich:Colimits}, hence so is the composite with \eqref{weil symmetric}.
The description \eqref{DXW} follows since the Kan extension is characterized by preserving limits in $\PreStk_\bk^\opp$.
\xpf

The limit in \eqref{DXW} can be equivalently computed in $\PrL$ or in $\Cat_\infty$. 
So, $\tD(X^W)$ is the categorical fixed point of $\varphi_X^*$ on $\D(\ol X)$.
More generally, for any finite family $(X_i)_{i\in I}$ in $\Sch_k^\ft$ the category $$\tD((X_i)_{i\in I}^W)=\tD\big(\prod_{i\in I}X_i^W\big)$$ is given by the simultanenous fixed points of $(\varphi_{X_i}^*)_{i\in I}$ on $\D(\prod_{i\in I}\ol X_i)$, cf.~\cite[Section~2.2]{HemoRicharzScholbach:Kuenneth}.

The categories $\tD(X^W)$ tend to be too large to admit a categorical Künneth formula.
For example this is the case for $\D$ being the derived category of étale torsion sheaves, and $X = Y = \A^1_k$, cf.~\cite[Remark~6.6]{HemoRicharzScholbach:Kuenneth}.
We therefore introduce the following subcategory.

\defi\thlabel{Weil formalism atomic}
Let $\D \colon (\Sch_\bk^\ft)^\opp \r \PrL$ be a lax symmetric monoidal functor that satisfies \thref{useful assumption}, i.e., takes values in atomically generated $\D(\bk)$-modules with functors preserving atomic objects.
For $X\in \Sch_k^\ft$, consider the full subcategory
\[
 \D(X^W)\subset \tD(X^W)
\]
that is the presentable full $\D(\bk)$-submodule generated by objects $M\in \tD(X^W)$ whose image $\ol M\in \D(\ol X)$ under the forgetful functor is $\D(\bk)$-atomic. 
Similarly, one defines $\D(\prod_{i\in I}X_i^W) \subset \tD(\prod_{i\in I}X_i^W)$ for a finite set $I$.
\xdefi

\rema
If $\D(\bk)$ is stable, compactly generated, and generated by dualizable objects, then being atomic is equivalent to being compact. 
Thus, $\D(X^W)$ is compactly generated by objects $(\ol M, \alpha: \ol M \stackrel \cong \r \varphi_X^* \ol M)$ with $\ol M \in \D(\ol X)^\omega$. 
The reader is referred to \cite[Definition~2.1 ff]{HemoRicharzScholbach:Kuenneth} for a further discussion of the difference between $\D(X^W)$ and $\tD(X^W)$.
\xrema

\lemm
\thlabel{DW lax tensor}

The assignment $(X_i)_{i\in I} \mapsto \D(\prod_{i\in I}X_i^W)$ in \thref{Weil formalism atomic} constitutes a lax symmetric monoidal functor
$$\D((-)^W)\colon (\Sch_k^\ft)^{\sqcap,\opp} \r \PrL.$$
In particular, for any $X,Y\in \Sch_k^\ft$, there is an exterior product functor
$$\D(X^W) \t_{\D(\bk)} \D(Y^W) \r \D(X^W\x_\bk Y^W). \label{schlax palax}$$
It is fully faithful if the exterior product functor for $\D$ on $\ol X,\ol Y$ is so.
\xlemm

\pf
By \thref{useful assumption}\eqref{useful assumption hom}, the exterior product functor for $\D$ preserves atomic objects.
Hence, $\DW(X)$ is a commutative $\D(\bk)$-subalgebra of $\tD(X^W)$.
As all pullback functors preserve atomic objects by assumption on $\D$ we obtain a lax symmetric monoidal functor as stated.

The full faithfulness of $\D(\ol X) \t_{\D(\bk)} \D(\ol Y) \r \D(\ol X\x_\bk \ol Y)$ implies the one for $\tD(X^W) \t_{\D(\bk)} \tD(Y^W) \r \tD(X^W\x_\bk Y^W)$, since mapping objects in limits of \ii-categories can be computed as limits of the mapping objects in the involved \ii-categories.
So, it suffices to see that $\D(X^W) \t_{\D(\bk)} \D(Y^W)$ is a full subcategory of $\tD(X^W) \t_{\D(\bk)} \tD(Y^W)$.
By definition, the atomic generators of $\D(X^W)$ are also atomic in $\tD(X^W)$, so we can apply \thref{internal left adjoint lemma} to the inclusion $\D(X^W) \subset \tD(X^W)$. 
Then, we apply \thref{proposition:infinity-cats-tensor-product} to see that $\D(X^W) \t_{\D(\bk)} \D(Y^W)$ is atomically generated and that it embeds fully faithfully into $\tD(X^W) \t_{\D(\bk)} \tD(Y^W)$.
\xpf

\lemm
\thlabel{D functor Weil}
Suppose that $\varphi_X^* : \D(\ol X) \r \D(\ol X)$ is an equivalence for all $X \in \Sch_k^\ft$.
Let $f \colon X \r Y$ be a map in $\Sch^\ft_k$ for which $\ol f^*\colon \D(\ol Y)\r \D(\ol X)$ is $\D(\bk)$-internally right adjoint, i.e., admits a $\D(\bk)$-linear left adjoint, denoted $\ol f_\sharp$. 
Then, the following holds:
\begin{itemize}
    \item $(\widetilde f^W)^* \colon \tD(Y^W) \r \tD(X^W)$ admits a $\D(\bk)$-linear left adjoint, denoted $(\widetilde f^W)_\sharp$.
    \item $(\widetilde f^W)_\sharp$ agrees on the level of the underlying categories with $\ol f_\sharp$.
    \item
    If $\D\colon (\Sch_\bk^\ft)^\opp\r \PrL$ satisfies \thref{useful assumption} and $\ol f_\sharp$ preserves $\D(\bk)$-atomic objects, then these functors restrict to an adjunction denoted by $(f^W)_\sharp : \DW(X) \rightleftarrows \DW(Y) : (f^W)^*$.
\end{itemize}
An identical statement holds for right adjoints, denoted $\ol f_*, (\widetilde f^W)_*$ and $(f^W)_*$ respectively, and also for maps in $(\Sch_k^\ft)^\sqcap$.
\xlemm

\pf
If $\ol f^*$ admits, say, a left adjoint, then the following diagram is vertically left adjointable since $\varphi_X^*$ is an equivalence:
$$\xymatrix{
\D(\ol Y) \ar[d]^{\ol f^*} \ar[r]^{\varphi_Y^*} & \D(\ol Y) \ar[d]^{\ol f^*} \\
\D(\ol X) \ar[r]^{\varphi_X^*} & \D(\ol X).
}$$
This implies that $\widetilde f_\sharp$ exists, commutes with $\ol f_\sharp$, and is $\D(\bk)$-linear.
If $\ol f_\sharp$ preserves atomic objects, then so does $\widetilde f_\sharp$, by definition of $\D((-)^W) \subset \tD((-)^W)$.
The same proof works if we replace ``right'' by ``left'' throughout.
The extension to $(\Sch_k^\ft)^\sqcap$ is straight forward. 
\xpf

\subsection{K\"unneth formula}\label{kuenneth formula weil}
Each of the formalisms $\D=\DA$, $\DAf$ or $\Dad$ on $\Sch_\bk^\ft$ takes values in atomically generated modules (\thref{DA compactly generated} and \thref{DM gm flat}), noting that $\bk$ is separably closed (\thref{example:cohomological-dimension}).
For $X \in \Sch_k^\ft$, the map $\varphi_X \colon \ol X \r \ol X$ is a universal homeomorphism, so that $\varphi_X^*$ is an equivalence on $\SH(\ol X)[p^{-1}]$ and on $\D(\ol X)$ for all three above-mentioned sheaf formalisms \cite[Theorem~2.1.1, Remark~2.2.9]{ElmantoKhan:Perfection}.
So, we can apply the discussion of the preceding section.
In the same vein as was discussed after \eqref{equation:definition-SH}, \thref{D functor Weil} and Zariski descent for $\DA$ allows to extend the construction in \thref{DW lax tensor} to a lax symmetric monoidal functor
$$\D((-)^W)\colon \Corr((\Sch_k^\ft)^\sqcap) \r \PrL$$
in each case $\D=\DA$, $\DAf$ or $\Dad$.
By construction, the functors $(f^W)^*, (f^W)_*$, $(f^W)_!, (f^W)^!$ agree with $\ol f^*, \ol f_*$, $\ol f_!, \ol f^!$ respectively after forgetting the Frobenius equivariance datum.

\lemm\thlabel{fully faithful Weil}
Let $\D=\DA$, $\DAf$ or $\Dad$.
Then, the exterior product functor 
\[
\D(X^W)\otimes_{\D(\bk)}\D(Y^W)\r \D(X^W\x_\bk Y^W)
\]
is fully faithful for any $X,Y\in \Sch_k^\ft$.
\xlemm
\pf
The exterior product functor for $\D$ and $\ol X, \ol Y$ is fully faithful by \thref{proposition:SH-fully-faithful} and \thref{fully faithful DAf}.
So, \thref{DW lax tensor} finishes the proof.
\xpf

The following lemma allows to leverage the categorical K\"unneth formula from ind-dualizable to all Weil objects.
Its analogue for $\ell$-adic Weil sheaves has been shown in \cite[Proposition~4.11]{HemoRicharzScholbach:Kuenneth}.

\prop
\thlabel{Zykellemma}
Let $\D = \DA$, $\DAf$ or $\Dad$.
For $X, Y \in \Sch_k^\ft$, the category $\D(X^W\x_\bk Y^W)$ is generated under colimits by the objects $(\iota_X\x \iota_Y)^W_! M$ for all locally closed embeddings $\iota_X \colon U \r X$, $\iota_Y \colon V \r Y$ and dualizable objects $M \in \D(U^W\x_\bk V^W)$.
\xprop
\pf
In each case, $\D(X^W\x_\bk Y^W)$ is atomically generated by objects $M$ whose underlying object $\ol M \in \D(\ol X \x \ol Y)$ is constructible: this follows from \thref{Weil formalism atomic} using \thref{characterization compact objects} for $\D=\DA$, \thref{atomic vs constructible} for $\D=\DAf$ and \thref{Dad comparison} for $\D=\Dad$.
Such objects $\ol M$ are dualizable when restricted to an open subset $T \subset \ol X \x \ol Y$.
The open subset $T' := \bigcup_{m, n \ge 0} (\varphi_{X}^{m}\x \varphi_{Y}^{n}) (T) \subset \ol X \x \ol Y$ is stable under the two partial Frobenii, and $\ol M|_{T'}$ is dualizable.
By \cite[Proposition~4.8]{HemoRicharzScholbach:Kuenneth} (this is the key geometric input), $T'$ contains an open subset of the form $\ol U \x \ol V$, for some open subsets $U \subset X$, $V\subset Y$.
Then, a Noetherian induction settles the claim.
\xpf

The following is the analogue of \thref{conjecture kuenneth} for étale Weil motives.
Interestingly, \thref{Zykellemma} allows to remove the additional assumptions present in the formulation of \thref{conjecture kuenneth}:

\conj\thlabel{conjecture kuenneth weil}
For any $X, Y \in \Sch_k^\ft$ the exterior product functor 
\eqn\label{conjecture kuenneth weil eq}
\DA(X^W) \t_{\DA(\bar k)} \DA(Y^W) \r \DA(X^W\x_\bk Y^W)
\xeqn
is an equivalence, and likewise on ind-dualizable categories $\D((-)^W)^\indbl$ (\thref{definition indbl}).
\xconj

The functor in \thref{conjecture kuenneth weil} is fully faithful by \thref{fully faithful Weil}. 
In addition, \thref{Zykellemma} implies that \eqref{conjecture kuenneth weil eq} is essentially surjective if it is so on ind-dualizable categories. 
The following result provides some evidence for the conjecture:

\theo
\thlabel{Künneth DMf Weil}
Let $\D = \DAf$ or $\Dad$.
For $X,Y\in \Sch_k^\ft$, the exterior product induces an equivalence
\begin{equation}
\D(X^W) \t_{\D(\bar k)} \D(Y^W) \stackrel \cong \r \D(X^W\x_\bk Y^W). \label{Künneth Weil}
\end{equation}
In other words, the induced functor
$$\D((-)^W) \colon \Corr((\Sch_k^\ft)^\sqcap) \r \Mod_{\D(\bk)}(\PrL)$$
is strictly (as opposed to lax) symmetric monoidal.
\xtheo

\rema
\thlabel{Zukunftsmusik} 
This statement can formally be extended to prestacks over $k$ (cf.~\thref{Künneth formula prestacks}).
\xrema

\pf
The functor \eqref{Künneth Weil} is fully faithful by \thref{fully faithful Weil}.
So, it is enough to show that the generators exhibited by \thref{Zykellemma} lie in the image of \eqref{Künneth Weil}.
Since $\bx$ is compatible with $!$-pushforward, it is enough to show that all objects $M\in \D(X^W\x_\bk Y^W)$ whose underlying object $\ol M\in \D(\ol X\x \ol Y)$ is dualizable lie in the essential image.
Furthermore, using localization (\thref{DMf localization}\eqref{DMf localization 1}) and a stratification into the essentially smooth loci of $X$, $Y$ and their complements, we may assume that both $X$ and $Y$ are smooth over $k$. 
In this case we claim that the functor
\[
\D(X^W)^\indbl \t_{\D(\bar k)} \D(Y^W)^\indbl \stackrel \cong \r \D(X^W\x_\bk Y^W)^\indbl
\]
is an equivalence, where $\indbl$ refers to the notation in \thref{definition indbl}.
The key difference to the non-Weil situation is that for $X$ smooth but \emph{not necessarily proper} over $k$, the pushforward along the projection $\pr_X : X \x Y \r X$ induces a functor 
$$(\pr_X^W)_* \colon \D(X^W\x_\bk Y^W) \r \D(X^W),$$
which preserves dualizable objects.
To see this note that $\pr_{X,*}$ preserves atomic objects, and that $M$ is dualizable if and only if $\ol M$ is dualizable.
So, if $\D=\DAf$ then we may check that claim by passing to $\Dad$: the adic realization $\rho^\flat$ is conservative, commutes with $\t$ and $\IHom$ between atomic objects, and therefore detects dualizability.
Further, as in the proof of \thref{pushforward char 0} we reduce to $\Z/n$-sheaves for some $n\in\Z_{\geq 1}$ prime to $p$.
By the categorical Künneth formula for $\Z/n$-sheaves proven in \cite{HemoRicharzScholbach:Kuenneth}, the dualizable object $M$ lies in the subcategory of $\Dad(X^W\x_\bk Y^W)/n$ generated under finite colimits and retracts by objects of the form $A \bx  B$, with $A \in \Dad(X^W)/n, B \in \Dad(Y^W)/n$ dualizable.
We compute $\pr_{X,*}(A \bx B) =  A \t p_{Y,*}B$, and this object is dualizable since any compact object in $\Dad(\bk)/n$ is dualizable.

Given this preservation of dualizability, and given that the $!$- and $*$-functors for $\D((-)^W)$ reduce to the classical functors for $\D$ on $\Sch_\bk^\ft$ we can repeat the proof of \thref{Künneth DMf indlis} in order to obtain the essential surjectivity of $\bx$ for $\Dad((-)^W)$ and $\DAf((-)^W)$.
\xpf

\bibliographystyle{alphaurl}
\bibliography{bib}

\end{document}